\documentclass{article}
\newcommand{\R}{\mathbb{R}}

\newcommand{\Rk}{R_k}
\newcommand{\Pk}{P_k}

\newcommand{\Pktilde}{\tilde{P}_k}

\newcommand{\Dktilde}{\tilde{D}_k}
\newcommand{\Vk}{V_k}
\newcommand{\Lambdak}{\Lambda_k}
\newcommand{\Zk}{Z_k}
\newcommand{\Tk}{T_k}
\newcommand{\E}{\mathbf{E}}
\newcommand{\Etilde}{\mathbf{\tilde{E}}}
\newcommand{\veps}{\varepsilon}
\newcommand{\VV}{\mathcal{V}}
\newcommand{\VVb}{\overline{\mathcal{V}}}
\newcommand{\Tktilde}{\tilde{T}_k}

\newcommand{\Mk}{\mathcal{M}_k}
\newcommand{\Mg}{\mathit{Mg}}
\newcommand{\Mgm}{\mathtt{Mg}}
\newcommand{\mA}{\mathtt{A}}
\newcommand{\mB}{\mathtt{B}}
\newcommand{\mC}{\mathtt{C}}
\newcommand{\mH}{\mathtt{H}}

\newcommand{\mV}{\mathtt{V}}
\newcommand{\mLa}{\mathtt{\Lambda}}
\newcommand{\ip}[1]{\langle{#1}\rangle}
\newcommand{\meas}{\mathrm{meas}}
\newcommand{\ev}{\mathtt{f}\!}
\newcommand{\cf}{\mathtt{e}}
\newcommand{\Hone}{H^1}
\newcommand{\Hmone}{H^{-1}}

\newcommand{\Hhalf}{H^{1/2}}
\newcommand{\Hmhalf}{H^{-1/2}}
\newcommand{\Hnkl}{H_0^{(1)}}

\usepackage{amssymb,amsthm,epsfig,amsmath}

\usepackage[dvips,colorlinks,bookmarksopen,bookmarksnumbered,citecolor=red,urlcolor=red]{hyperref}
\newtheorem{theorem}{Theorem}
\newtheorem{lemma}{Lemma}
\newtheorem{algorithm}{Algorithm}

\begin{document}

\title{Convergence analysis of a multigrid algorithm for the acoustic single layer equation}

\author{S.Gemmrich\footnote{Dept. of Mathematics and Statistics, McGill University}, J.Gopalakrishnan\footnote{Dept. of Mathematics, Portland State University}, N.Nigam\footnote{Dept. of Mathematics and Statistics, McGill University}}
\maketitle
\begin{abstract}

We present and analyze a multigrid algorithm for the acoustic single layer equation in two dimensions. The boundary element formulation of the equation is based on piecewise constant test functions and we make use of a weak inner product in the multigrid scheme as proposed in \cite{BLP94}. A full error analysis of the algorithm is presented. We also conduct a numerical study of the effect of the weak inner product on the oscillatory behavior of the eigenfunctions for the Laplace single layer operator.

This paper is dedicated to Prof. G.C.Hsiao on the occasion of his 75th birthday.
\end{abstract}

\section{Introduction}\label{nlg_intro}
A model for the scattering of acoustic waves by a bounded obstacle is given by the Helmholtz equation in the exterior of the scatterer, with appropriate growth conditions on the scattered field. One can reformulate this problem in terms of integral equations on the surface of the scattering object via direct or indirect boundary integral formulations. 
Or one may consider scattering from a coated bounded obstacle, in which case an integral equation can be used to prescribe a non-reflecting condition on an artificial surface surrounding the object. In both, one has to find numerical approximations to solutions of boundary integral equations. In this context, Galerkin type methods have been studied extensively and become popular over recent years, see for example the monographs \cite{SS04}, \cite{HW08}.

Our focus in this paper lies on integral equations of the first kind, which arise naturally in the direct  boundary integral method for the Dirichlet problem.  The main integral operator involved is called the single layer operator and may be viewed as a pseudo-differential operator of order minus one. Several authors have observed advantages of using integral equations of the first kind (e.g. \cite{HM73}), for example when the scattering object is very thin (e.g. \cite{HSW91}) or when the scattering surface is not smooth. Indeed, for problems of crack propagation in elasticity or scattering from a screen, integral equations of the first kind are the most appropriate model. 

Other popular integral equation strategies include the use of combinations of single and double layer operators, to avoid issues of non-uniqueness and the potential for numerical instability near possible eigenvalues of the operators. These approaches include the famous Brakhage-Werner and Burton-Miller formulations.  In these cases as well, discrete strategies rely on effective and accurate approximation methods for the layer operators involved.

Due to the non-local behavior of  boundary integral operators, they typically lead to dense linear systems upon discretization. Though one only needs to mesh on a surface of co-dimension one, the fill-in in the matrices corresponding to the integral operators is significant. Without some form of preconditioning or acceleration, these methods then become prohibitive. 

One possible preconditioning strategy is the use of a multigrid scheme. However, the use of standard multigrid smoothing operations is inappropriate for negative order pseudodifferential operators. Such operators link highly oscillatory eigenfunctions to small magnitude eigenvalues. This ruins the basic interplay between standard smoothing of oscillatory error components and the possibility to represent the remaining error components on coarser grids. The remedy for this is the use of a weaker inner product in order to modify this spectral feature of the operator. This approach has first been described in \cite{BLP94} for positive-definite operators. We follow the same path for the acoustic single layer equation. The numerical examples in Section~\ref{nlg_sec::ne} exemplify how the spectral behavior of the discretized operators changes through the use of the weaker inner products.

Several related works can be found in the literature. The single layer equation associated to Laplace's equation has been treated and analyzed in \cite{FS97} using a BPX preconditioner. The same equation was considered in \cite{LP07}. Here, the authors studied a multigrid method for large-scale data-sparse approximations to the single layer operator.
In the acoustics case, the use of Haar basis functions and compression type multilevel algorithms for integral equations of the first kind has been studied in \cite{MMS97}. Algebraic multigrid preconditioners, based on the smoother in \cite{BLP94}, have been developed in \cite{LPR03}. 

The purpose of this paper is to prove convergence of the multigrid algorithm given in \cite{BLP94} when applied to the indefinite acoustic single layer discretization. We will make use of perturbation-type arguments such as 
in \cite{BKP94} and \cite{GPD04}. We also state the algorithm in the situation where a non-uniform discretization is used.

The design of the algorithm heavily relies on the above references and we have reported its promising numerical 
performance in \cite{GGN08}.
Currently, our codes assemble and multiply  matrices in 
O($n^2$) complexity, and so the computational cost of the iterative solution by multigrid is the 
same. The efficient assembly and matrix multiplication of the corresponding matrices is an active research issue in multipole and 
hierarchical matrix theories. Our numerical results indicate that we can approach optimality in complexity 
through multigrid, once assembly and matrix multiplications are done 
optimally.

At this juncture, we note important directions for extending this work, motivated by recent developments in the study of boundary integral equation methods. Firstly, the use of combined integral equation formulations, including those by Brakhage-Werner \cite{BrakhageWerner} and Burton-Miller, provide stable formulations for the solution of scattering problems. These are popular approaches, and the development of a multigrid strategy for these would be of interest. Both the methods of implementation and analysis would be different than those in this paper, to account for the different properties of the combined layer operators. Secondly, the methods presented in this paper are not tailored to high-frequency scattering, and our analysis does not include wave-number explicit bounds. Wave-number explicit bounds are described, for example, for  combined field operator approaches \cite{snc01,snc02,snc03,melenk2010, melenk2011}. The scale resolution condition of $kh$ being sufficiently small (for piecewise linears) needs to be met;  in high-frequency settings non-polynomial approximation spaces may be a better approach, \cite{ansatz}.  The dependence of the coarsest mesh on the wave number is assumed to satisfy this requirement. In this paper, our focus is on a much simpler situation:  how to use a first-kind integral equation with piecewise constant approximants to the solution, in the low to medium frequency situations of scattering from polygonal domains.

We now give a brief derivation of the acoustic single layer equation using the framework of a direct boundary integral approach. We consider the following exterior Helmholtz problem with prescribed Dirichlet data on the boundary of a scatterer. Here, $\Gamma$ is a simple, closed polygonal curve in the plane and  $\Omega^{ext}$ denotes its exterior domain. 
\begin{align*} 
	-\Delta u - \kappa^2 u	 =  0 	\quad \mbox{ in $\Omega^{ext}$},\qquad u  =  g \quad \mbox{ on $\Gamma$}\qquad\mbox{and} \qquad
 \lim_{r\to\infty} r^{\frac{1}{2}} (\frac{\partial u}{\partial r}-i\kappa u)  =  0.
\end{align*}

To guarantee unique solvability we assume a non-zero, real wave-number $\kappa\in\mathbb{R}$, such that $\kappa^2$ is not an interior eigenvalue of $-\Delta$. 
The Sommerfeld radiation condition is given in terms of the usual radial component $r$ in polar coordinates. 
It is well known that the solution to this problem is fully determined by its complete Cauchy data $g=\gamma^{+}u$ and $\sigma=B_{\nu}^{+}u$, where $\gamma^{+}: H^1_{\textup{loc}}(\Omega^{ext})\rightarrow H^{1/2}(\Gamma)$ and $B_{\nu}^{+}: H^1_{\textup{loc}}(\Omega^{ext})\rightarrow H^{-1/2}(\Gamma)$ denote the exterior trace operator and the exterior conormal derivative respectively. The normal $n$ is assumed to be outward to $\Omega^{ext}$, ie, it points into the bounded region.
In fact, for $x \in \Omega^{ext}$ one has an integral representation formula for the solution to the boundary value problem (e.g.\ \cite{HW08}, \cite{S03}), namely
\begin{align}\label{nlg_eq:repr}
u(x) & = -\frac{i}{4} \int_{\Gamma}\Hnkl(\kappa |x-y|)\, \sigma(y)  ds_y\,+\, \frac{i}{4}\int_{\Gamma}\frac{\partial \Hnkl(\kappa |x-y|)}{\partial n_y} \, g(y)  ds_y .
\end{align}  
The kernels of the two integrals are given in terms of the Hankel function $\Hnkl(z)$ and its conormal derivative.
According to the representation~\eqref{nlg_eq:repr}, it is sufficient to find the unkown surface density $\sigma$. To do this, one exploits the jump relations of the two integrals, i.e.\ their behavior in the limit as $x$ approaches $\Gamma$ from both the interior and the exterior of the scattering domain. These relations appear when we take the trace of equation \eqref{nlg_eq:repr} and lead to the following integral equation for $\sigma$:
\begin{align}
\label{nlg_eq:sl}
V \, \sigma &= f \in\Hhalf(\Gamma).
\end{align}
Here, the right hand side $f=(\frac{1}{2}I+K)\, g$ depends on the Dirichlet trace $g$ and requires the evaluation of the double layer operator $K$. The single layer operator $V$  and the double layer operator $K$ are both defined in terms of singular integrals.
\begin{align*}
 &V: \Hmhalf(\Gamma) \rightarrow \Hhalf(\Gamma), & V\sigma(x) &:= \frac{i}{4} \int_{\Gamma}\Hnkl(\kappa |x-y|)\, \sigma (y) ds_y  & x\in\Gamma, \\
 &K:\Hhalf(\Gamma) \rightarrow \Hhalf(\Gamma), & K\mu(x) &:= \frac{i}{4} \int_{\Gamma}\frac{\partial \Hnkl(\kappa |x-y|)}{\partial n_y}\,  \mu (y) ds_y  & x\in\Gamma.
\end{align*}

One should note that the above approach to reformulate the original boundary value problem into an integral equation on $\Gamma$ is by no means unique. However, several methods lead to a single layer equation of the form \eqref{nlg_eq:sl}.
In the following, we are interested in the weak form of equation \eqref{nlg_eq:sl}.
Given $f\in \Hhalf(\Gamma)$, find $\sigma\in \Hmhalf(\Gamma)$ such that
\begin{align}
\label{nlg_eq:weak}
\mathcal{V}(\sigma,\mu) &=\langle f, \mu\rangle \quad\mbox{ for all
}\mu\in
\Hmhalf(\Gamma),
\end{align}
where the continuous sesquilinear form $\mathcal{V}
:H^{-1/2}(\Gamma)\times H^{-1/2}(\Gamma)\rightarrow \mathbb{C}$ is  defined
by 
$\mathcal{V}(\sigma,\mu)= \langle V\sigma,\mu\rangle$,
and $\langle{\cdot,\cdot}\rangle$ denotes
the duality pairing between
$\Hhalf(\Gamma)$ and $\Hmhalf(\Gamma)$.  
The single layer operator which corresponds to the Laplacian will be denoted by $\Lambda$.
\begin{align*}
&\Lambda: \Hmhalf(\Gamma) \longrightarrow \Hhalf(\Gamma), & \Lambda\sigma(x) &:= -\frac{1}{2\pi} \int_{\Gamma} \ln(|x-y|)\, \sigma (y) ds_y  & x\in\Gamma.
 \end{align*}
We note that the underlying differential operator is the principal part of the Helmholtz operator.
Its associated sesquilinear form $\Lambda(\cdot,\cdot)$, defined similarly to $\mathcal{V}(\cdot,\cdot)$ above, is positive definite after the region has been scaled properly (see \cite{S03} for details), i.e. there holds
\begin{align}\label{nlg_eq:posdef}
 \Lambda (\sigma ,\sigma) &\geq C\, \|\sigma\|_{\Hmhalf(\Gamma)}^2 \qquad\mbox{ for all } \sigma\in\Hmhalf(\Gamma).
\end{align}
Consequently, it defines an inner product whose induced norm is equivalent to the natural energy norm in ${\Hmhalf(\Gamma)}$.
This will play an important role in the analysis of Section~\ref{nlg_sec:conv}.

The paper is organized as follows: In Section~\ref{nlg_sec::ma} we present a multigrid algorithm for integral equations of the first kind, and introduce a (computable) inner product whose induced norm is equivalent to the natural norm in $\Hmone(\Gamma)$ on finite dimensional test spaces. The multigrid strategy relies on reformulating Problem~\eqref{nlg_eq:weak} using both the standard inner product in $\Hmone(\Gamma)$ and the new computable version for piecewise constant functions. We introduce smoothers, and present a matrix version of the algorithm. Section~\ref{nlg_sec:conv} consists of a convergence analysis. The key component is a careful study of the difference between the single layer operators for the Laplacian and for the Helmholtz equation. We conclude this section with a convergence result. Finally, in Section~\ref{nlg_sec::ne} we present some numerical experiments describing the spectral behaviour of the single layer operator, as well as that of the operators used in the reformulated problem. We see, in the context of a smooth curve (a circle) and a Lipschitz curve (a square) how the use of the weaker inner product renders the problem suitable for a multigrid strategy. We conclude by reporting on the convergence behavior of the algorithm applied to simple test cases.

%%%%%%%%%%%%%%%%%%%%%%%%%%%%%%%%%%%%%%%%%%%%%%%%%%%%%%%%%%%%%%
%%%%%%%%%%%%%%%%%%%%%%%%%%%%%%%%%%%%%%%%%%%%%%%%%%%%%%%%%%%%%%
%%    a multigrid algorithm   
%%%%%%%%%%%%%%%%%%%%%%%%%%%%%%%%%%%%%%%%%%%%%%%%%%%%%%%%%%%%%%
\section{A multigrid algorithm}
\label{nlg_sec::ma}

In this section we present the multigrid algorithm, originally proposed
in~\cite{BLP94} for the positive definite pseudodifferential operators
of order minus one, and which was applied to the acoustic case
in~\cite{GGN08}.  The multigrid algorithm is presented below in
terms of a smoother, whose definition is postponed to
Subsection~\ref{nlg_ssec:smooth}.  The smoother is realized using a weak
base inner product.  We make this more precise in
Subsection~\ref{nlg_sec:discr}.

We first establish notations needed to describe the multigrid
algorithm. Assume that the polygonal boundary $\Gamma$ is composed of
finitely many straight edges $\Gamma_j$. Each $\Gamma_j$ is meshed by
a coarse grid of line segments of length $l_j^1$.

We successively refine this grid in a uniform way by breaking each
element in half and adding the respective midpoints to the vertices of
the next finer level.  On every level of refinement $k=1,2,\ldots ,
J,$ we label the vertices in such a way that $x^k_1,x^k_2,\ldots
,x^k_{N_k},x^k_{N_k+1}=x_1^k$ is a counterclockwise enumeration. Now,
let $\phi^k_i$ be the characteristic function of the line segment
$\tau^k_i=\mathrm{conv}(x^k_i,x^k_{i+1})$ ($i=1,2,\ldots ,N_k$) and
denote their span by $\Mk:=\mathrm{span}\left\{\phi^k_i\right\}$. For
the sake of easier notation we will suppress the level number $k$ in
our notations whenever the context rules out any ambiguity.  This
construction yields a sequence of nested finite-dimensional spaces
\[ \mathcal{M}_1 \subset \mathcal{M}_2 \subset\ldots\subset\mathcal{M}_J \subset \Hmhalf(\Gamma).\]

We now define the discrete operators $V_k:\Mk\rightarrow\Mk$ with the
help of the $\Hmone (\Gamma)$ inner product, denoted by
$(\cdot,\cdot)_{-1}$. The defining relation is
\begin{align}
\label{nlg_eq:1}
(V_k \sigma,\mu)_{-1} & =\mathcal{V}(\sigma,\mu) \quad \mbox{ for all } \sigma, \mu\in\Mk.
\end{align}
Analogously, we choose $f_k\in\Mk$ to satisfy
$(f_k,\mu)_{-1} =\,\langle f , \mu \rangle$  
for all $ \mu\in\Mk.$
Then, on every level $k$, the equation of interest can be written in
operator form as 
\begin{align}\label{nlg_eq:2}
V_k \,\sigma_k & = f_k.  
\end{align}
In order to describe the algorithm in a function space setting, we shall also need the $\Hmone$-
projections $Q_k:\Hmone(\Gamma)\rightarrow\Mk$, which are defined by 
\begin{align*}
(Q_k \sigma,\mu)_{-1} &=(\sigma,\mu)_{-1} \quad \mbox{for all}\; \mu\in\Mk.
\end{align*}
We further need a family of smoothing operators
$R_k:\Mk\rightarrow\Mk$. It is defined precisely in
Subsection~\ref{nlg_ssec:smooth}, but for now, we just assume that $R_k$
are some given linear operators.  Then, given an initial guess
$\sigma_0\in\mathcal{M}_J$, the multigrid iteration computes a
sequence of approximate solutions to \eqref{nlg_eq:2} using an iteration
of the form $ \sigma_{i+1}=\Mg_J(\sigma_i,f_J )$, where $\Mg_J(\cdot\,
, \cdot)$ is a mapping of $\mathcal{M}_J\times\mathcal{M}_J$ into
$\mathcal{M}_J$, defined recursively by the following algorithm:

\begin{algorithm}\label{nlg_a:1}
 Set $\Mg_1(\sigma,f)={V_1}^{-1}f$. If $k>1$ we define
 $\Mg_{k}(\sigma, f)$ recursively as follows:
\begin{align}
\label{nlg_eq:13}
\sigma_1		& =\sigma	 + R_k(f-V_k\sigma), \\
\label{nlg_eq:14}
 \Mg_k(\sigma,f)& =\sigma_1 + \Mg_{k-1}(0, Q_{k-1} (f-V_k\sigma_1)).
\end{align}
\end{algorithm}
This is a simple variant of a V-cycle multigrid scheme, which only
uses pre-smoothing. Equivalently, we can write the iterative scheme as a linear iteration method
\[ \sigma_{i+1}=\sigma_i\,+\,B_J\,(f_J-V_J\,\sigma_i ),\]
with an ``approximate inverse'' $B_J: \mathcal{M}_J \mapsto \mathcal{M}_J$ defined by
\[
B_k f_k = \mathit{Mg}_k(0,f_k) \quad \text{ for all } f_k \in \mathcal{M}_k
\text{ and } k=2,\ldots, J.
\]
This operator is useful as a preconditioner in preconditioned
iterative methods. The matrix version of Algorithm~\ref{nlg_a:1} is given
in Subsection~\ref{nlg_ssec:matrix}.

%%%%%%%%%%%%%%%%%%%%%%%%%%%%%%%%%%%%%%%%%%%%%%%%%%%%%%%%%%%%%%
%%    discrete inner products   
%%%%%%%%%%%%%%%%%%%%%%%%%%%%%%%%%%%%%%%%%%%%%%%%%%%%%%%%%%%%%%
\subsection{Discrete inner products}\label{nlg_sec:discr}

The use of the $\Hmone(\Gamma)$ inner product in defining the
operators $V_k$ and $Q_k$ for the multigrid algorithm confronts us with the question of
computability. We will have to work around this issue by introducing
equivalent computable inner products. These inner products will be used
to define smoothers in Section~\ref{nlg_ssec:smooth}.

The problem of calculating the $\Hmone(\Gamma)$ inner product of two elements
in $\Mk$ is related to the solution operator of a second order
boundary value problem on the boundary curve, namely
\begin{align}\label{nlg_MinusOne}
-u''+u=v.
\end{align} 
Here, 
$v\in\Hmone(\Gamma)$ is a given function, 
$u$ has  periodic boundary conditions,
 and the primes denote
differentiation with respect to arc-length. The weak form of this
problem is uniquely solvable and we denote the bounded solution
operator by $T:\Hmone(\Gamma) \longrightarrow \Hone(\Gamma)$.  Then,
for $v,w\in \Hmone(\Gamma)$ it is easily verified that
$(v,w)_{-1}=(Tv,w)_{\Gamma}=(v,Tw)_{\Gamma}$, where
$(\cdot,\cdot)_{\Gamma}$ denotes the complex $L^2(\Gamma)$-inner product.
Unfortunately, the use of the exact solution operator $T$ is
infeasible. Instead, we discretize \eqref{nlg_MinusOne} using a
second-order finite difference method.  This finite difference method
results in an $N_k\times N_k$ linear system with an tridiagonal-like
matrix $\mA_k$.  We need to introduce some more notation in order to
see the details. Functions in $\Mk$ can be represented through their
basis expansions with respect to the $\phi_j^k$. This is done via the
following map:
\begin{align}
  \label{nlg_eq:5}
  \mathtt{e}_k
  &: \Mk\longrightarrow \mathbb{C}^{N_k}, 
  && [\mathtt{e}_k (\sigma)]_i=  \frac{1}{l_i^k}\;(\sigma,\phi^k_i)_\Gamma, 
\end{align}
where $l_i^k:=\meas({\tau^k_i})$. Since the basis functions $\phi^k_i$
are the (orthogonal) indicator functions of the segments $\tau^k_i$,
the basis expansion for any $\sigma\in\Mk$ is $
\sigma=\sum_{l=1}^{N_k} \,[\mathtt{e}_k(\sigma)]_l\, \phi_l^{k},$ so
the map in \eqref{nlg_eq:5} gives the vector of coefficients.

We then define the (invertible) operator $A_k: \Mk\longrightarrow\Mk$ through
\begin{align}\label{nlg_eq:defA}
 A_k \sigma &= \sigma \;-\;\sum_{i=1}^{N_k}\left(\frac{[\cf(\sigma)]_{i+}-[\cf(\sigma)]_{i}}{l_i^2} \;-\; \frac{[\cf(\sigma)]_{i}-[\cf(\sigma)]_{i-}}{l_i\,l_{i-}}\right)\;\phi_i^k ,
\end{align}
where we use the notation $i+=\left(i+1\mod N_k\right)$ and
$i-=\left(i-1\mod N_k\right)$. As in~(\ref{nlg_eq:defA}), we will often
drop the superscript $k$ and identify $l_i^k$ and $l_i$ to be the
same for convenience.

The inverse operator $A_k^{-1}$ or equivalently the inverse matrix
$\mathtt{A_k^{-1}}$ serve as approximations for the solution
operator~$T$.  This motivates the definition of the computable,
discrete inner products on the spaces $\Mk$ via
\begin{equation}
\label{nlg_eq:12}
[\phi,\psi]_{k}:=\left(A_k^{-1} \phi,\psi\right)_{\Gamma}
\quad\mbox{for all } \;\phi ,\psi \in\mathcal{M}_{k}.
\end{equation}

The following lemma shows that $\|\cdot\|_{-1}$ and $[\cdot,\cdot]_k$
are equivalent norms, with the equivalence constants independent of
the refinement levels $k=1,..., J$. This result can be found in
\cite{BLP94} and \cite{FS97} for the case of the screen problem when
$\Gamma$ is in fact an open boundary patch or an open line segment respectively. Here, we give a proof for a
closed boundary curve $\Gamma$.  For the analysis we will assume that
all the mesh element lengths $l^k_i$ are such that $C_1 h_k \le l^k_i
\le C_2 h_k$ for some fixed positive constants $C_1,C_2.$ Here $h_k$
is a representative mesh size, e.g., $h_k = \max( \meas(\tau_i^k))$.

For the proof we  need a standard approximation result. Let
$\theta\in\Hone(\Gamma)$ and let $\theta_k\in\Mk$ denote the piecewise
constant approximation defined by
\begin{align}
\label{nlg_eq:6}
 \theta_k &:= \sum_{i=}^{N_k} \theta(x^k_i) \,\phi_i^k.
\end{align}
This is well defined for $\theta$ in $\Hone(\Gamma)$ by the Sobolev
inequality
\begin{equation}   \label{nlg_eq:Sobolev}
 \|\theta\|_{L_{\infty}(\Gamma)}\leq C\, \|\theta\|_{\Hone(\Gamma)}  \quad\text{ for all }\theta
 \in \Hone(\Gamma).
\end{equation}
Then
\begin{equation}
  \label{nlg_eq:approx}
  \| \theta - \theta_k \|_{L^2(\Gamma)}\le C h_k |\theta
  |_{\Hone(\Gamma)}
  \quad\text{ for all }\theta
 \in \Hone(\Gamma).
\end{equation}
To prove this, it suffices to observe that 
\begin{align*}
 \|\theta -\theta_k \|^2_{L^2(\tau_i^k)} &= \int_{\tau_i^k} |\theta
 (x) - \theta(x^k_i)|^2 \; dx
 = \int_{\tau_i^k} \left| \int_{x^k_i}^x \theta^{\prime}(\xi)\,d\xi \right|^2 dx\\
 &
 \leq \int_{\tau_i^k} \left( \int_{x^k_i}^{x^k_{i+1}} |\theta^{\prime}(\xi)|^2\,d\xi\right) \left(x-x^k_i \right) dx\\
        & = \frac{l_i^2}{2} \left(\int_{\tau_i^k} |\theta^{\prime}(\xi)|^2\, d\xi \right) = \frac{l_i^2}{2}\, |\theta|_{\Hone(\tau_i^k)}^2.
\end{align*}
Equation~(\ref{nlg_eq:approx}) follows by summing over all elements. We
will use~(\ref{nlg_eq:approx}) in the proof of the next lemma.

\begin{lemma}\label{nlg_lemma:contrib}
There exist constants $c,C>0$ independent of $J$ such that
\begin{align}\label{nlg_eq:innerpr}
c\|\sigma\|_{-1}^2 &\leq [\sigma ,\sigma ]_k \leq C\|\sigma\|_{-1}^2\qquad\mbox{for}\quad \sigma\in\Mk \quad (k=1,\ldots, J).
\end{align}
\end{lemma}
\proof{ Given an element 
\[
\sigma =\sum_{i=1}^{N_k}[\cf(\sigma)]_i  \,\phi^k_i,
\quad\text{ we define }\quad
\tilde\sigma=\sum_{i=1}^{N_k}
[\cf(\sigma)]_i
  \,\psi^k_i,
\]
where $\psi^k_i$ is the continuous and piecewise linear function which
takes the value one at node $x^k_i$ and vanishes at every other node
(otherwise known as the ``hat'' function). Note that $\tilde \sigma$
is in $\Hone(\Gamma)$. Then
applying~(\ref{nlg_eq:approx}) to $\tilde\sigma$, we obtain
\begin{align*}
\|\sigma - \tilde\sigma\|_{L^2(\Gamma)}\leq C\, h_k\,\|\partial \tilde\sigma\|_{L^2(\Gamma)},
\end{align*}
where $\partial\tilde\sigma$ denotes the derivative of $\tilde\sigma$
with respect to arc length.  Note that
\begin{align}
\label{nlg_eq:22}
\partial \tilde\sigma = 
\sum_i \left( \frac{[\cf(\sigma)]_{i+}-[\cf(\sigma)]_i}{l_{i}}\right)\phi^k_{i}.
\end{align}
From this it follows by straight forward calculations that
\begin{align}\label{nlg_eq:eingefuegt}
| \tilde\sigma|_{H^1(\Gamma)}^2 
% \leq (1 + C \,h_k^{-2})\,\|\sigma \|_{L^2(\Gamma)}^2 
\le C\, h_k^{-2}\, \| \sigma\|_{L^2(\Gamma)}^2.
\end{align}

Before proving the inequalities, we make a few observations. Using the
operator from \eqref{nlg_eq:defA} we see for $\sigma ,\mu\in \Mk$,
\begin{align*}
&(A_k \sigma , \mu)_{\Gamma} \\
  & = (\sigma ,\mu)_{\Gamma}
  - \sum_{i}\int_{\tau_i} \left(
  \frac{[\cf(\sigma)]_{i+}-[\cf(\sigma)]_i}{l_i^2}\phi_i -
  \frac{[\cf(\sigma)]_{i}-[\cf(\sigma)]_{i-}}{l_i l_{i-}}\phi_i\right) 
       \overline{[\cf(\mu)]}_i \phi_i\; ds\\
&		
=(\sigma ,\mu)_{\Gamma}
 - 
\sum_{i} \left( \frac{[\cf(\sigma)]_{i+}-[\cf(\sigma)]_i}{l_i} -  							\frac{[\cf(\sigma)]_{i}-[\cf(\sigma)]_{i-}}{ l_{i-}}\right) \overline{[\cf(\mu)]}_i \\
&
=(\sigma ,\mu)_{\Gamma} 
 + \sum_{i} \left( [\cf(\sigma)]_{i+} - [\cf(\sigma)]_i \right)\left(  											\overline{[\cf(\mu)]}_{i+}-\overline{[\cf(\mu)]}_{i}\right)  \frac{1}{l_i}\\
&		=(\sigma ,\mu)_{\Gamma} + \int_{\Gamma}
       \partial\tilde\sigma \,
       \overline{\partial\tilde\mu} \, ds.
\end{align*}
by~(\ref{nlg_eq:22}). The right hand side above defines a sesquilinear form 
\[
 a(\tilde\sigma ,\tilde\mu)
= (\sigma ,\mu)_{\Gamma} + 
(\partial\tilde\sigma  ,  \partial\tilde\mu)_{\Gamma}.
\]
It is easy to see that 
$\|\sigma\|^2_{L^2(\Gamma)} \sim \|\tilde\sigma\|^2_{L^2(\Gamma)}$, 
where $X\sim Y$ indicates that $ C_1 X \le Y \le C_2 Y$ holds with
some positive constants $C_1$ and $C_2$ independent of the mesh.
Thus we have proven that
\[
\| \tilde \sigma \|^2_{\Hone(\Gamma)} \,\sim \;
a(\tilde\sigma,\tilde\sigma) \equiv
(A_k\sigma,\sigma), 
\quad \text{for all } \sigma \in \Mk.
\]
If we combine the above norm equivalence with \eqref{nlg_eq:eingefuegt} and write $\|\sigma \|_{A_k}^2  = (A_k\sigma,\sigma)$, it immediately follows that 
\begin{align*}
 \|\sigma \|_{A_k}^2 
 &
 = a(\tilde\sigma,\tilde\sigma)  
 = \| \sigma \|_{L^2(\Gamma)}^2 +  | \tilde\sigma|_{H^1(\Gamma)}^2 
 \leq (1 + C \,h_k^{-2})\,\|\sigma \|_{L^2(\Gamma)}^2 
\le C\, h_k^{-2}\, \| \sigma\|_{L^2(\Gamma)}^2.
\end{align*}
Here, $C$ denotes a generic constant independent of meshsize. 
In other words, this yields the inequality 
\begin{equation}
  \label{nlg_eq:max}
  \lambda_{\mathrm{max}}(A_k) \le C h_k^{-2},
\end{equation}
which could also be shown using other methods.

We now begin proving  the inequalities of the lemma, starting
with the second inequality in \eqref{nlg_eq:innerpr}.
\begin{align}
\nonumber
(A_k \sigma , \sigma)_{\Gamma} &
=a(\tilde\sigma ,\tilde\sigma)= \sup_{\mu\in\Mk} \frac{|a(\tilde\sigma
  ,\tilde\mu)|^2}{a(\tilde\mu ,\tilde\mu)} \leq C \sup_{\mu\in\Mk}
\frac{|(A_k \sigma , \mu)_{\Gamma}|^2}{C_1
  \|\tilde\mu\|^2_{\Hone(\Gamma)}}\\
\nonumber
&  
\leq C \sup_{\mu\in\Mk} 
\frac{|(A_k \sigma , \tilde\mu)_{\Gamma}|^2+|(A_k \sigma ,
  \mu-\tilde\mu)_{\Gamma}|^2}{\|\tilde\mu\|^2_{\Hone(\Gamma)}}
\\ \nonumber
&  \leq C \sup_{\mu\in\Mk} \frac{|(A_k \sigma ,
  \tilde\mu)_{\Gamma}|^2}
  {\|\tilde\mu\|^2_{\Hone(\Gamma)}}+
  \frac{\|A_k \sigma\|^2_{L^2(\Gamma)} \|
    \mu-\tilde\mu\|^2_{L^2(\Gamma)}}
       {\|\tilde\mu\|^2_{\Hone(\Gamma)}}\\ \nonumber
	&\leq C \sup_{\mu\in\Hone(\Gamma)} \frac{|(A_k \sigma ,
         \mu)_{\Gamma}|^2}
       {\|\mu\|^2_{\Hone(\Gamma)}}\,+\, 
       C \,h_k^2\, \|A_k \sigma\|^2_{L^2(\Gamma)}\\
\label{nlg_eq:8}
	&= C\left(\|A_k \sigma\|^2_{\Hmone(\Gamma)} \,+\,h_k^2\, \|A_k \sigma\|^2_{L^2(\Gamma)}\right).
\end{align}
The inverse inequality 
\begin{equation}
  \label{nlg_eq:inv}
  \| \mu \|_{L^2(\Gamma)} \le C h_k^{-1} \| \mu \|_{\Hmone(\Gamma)}
  \quad\text{ for all } \mu \in \Mk
\end{equation}
can be found in~\cite[Eq.~(3.20)]{BLP94}. Using it in~(\ref{nlg_eq:8}),
we get the upper bound
\begin{align}
(A_k \sigma , \sigma)_{\Gamma}&\leq C \, \|A_k \sigma\|^2_{\Hmone(\Gamma)} \qquad\mbox{ for all $\sigma\in\Mk$}
\intertext{or equivalently}
(A_k^{-1} \mu , \mu)_{\Gamma}&\leq C \, \| \mu\|^2_{\Hmone(\Gamma)} \qquad\mbox{ for all $\mu\in\Mk$}.
\end{align}
This is the upper inequality stated in the lemma. 

It remains to prove the first inequality in~\eqref{nlg_eq:innerpr}. For
this, we need a stability result for the $L^2$-orthogonal
projection $\varPi_k:\Hone(\Gamma)\rightarrow\Mk$, namely
%% We  note that
%% \begin{align*}
%% (A_k^{-1} \sigma, \sigma ) &= \sup_{\mu\in\Mk}\frac{(\sigma , \mu)^2}{(A_k \mu ,\mu)}.
%% \end{align*}
\begin{align}\label{nlg_eq:QkProj}
 (A_k \varPi_k \theta ,\varPi_k\theta)&\leq C\, \|\theta\|_{\Hone(\Gamma)}^2.
\end{align}
The result \eqref{nlg_eq:QkProj} follows once we prove that for all
$\theta\in\Hone(\Gamma)$, there exists $\theta_k\in\Mk$, such that
\begin{align}
 \| \theta_k \|_{A_k} &\leq C\, \|\theta\|_{\Hone(\Gamma)} \label{nlg_eq:claim2i} \\
 \|\theta -\theta_k \|_{L^2(\Gamma)} &\leq C\, h_k\, |\theta|_{\Hone(\Gamma)}. \label{nlg_eq:claim2ii}
\end{align}
Indeed, writing $\varPi_k\theta =\varPi_k (\theta-\theta_k) +\theta_k$,
\begin{align*}
 \|\varPi_k\theta\|_{A_k} &\leq \|\varPi_k (\theta -\theta_k )\|_{A_k} +  \|\theta_k\|_{A_k}\\
        & \leq C\, h_k^{-1}\, \|\varPi_k (\theta -\theta_k )\|_{L^2(\Gamma)} +
 \|\theta_k\|_{A_k}  && \text{ by~(\ref{nlg_eq:max})}\\
        & \leq C \, h_k^{-1}\,\|(\theta -\theta_k )\|_{L^2(\Gamma)} + C\,
 \|\theta\|_{\Hone(\Gamma)}
 && \text{ by~(\ref{nlg_eq:claim2i})}\\
        &\leq C \, \|\theta\|_{\Hone(\Gamma)} &&\text{ by~(\ref{nlg_eq:claim2ii})}.
\end{align*}
Therefore, let us now exhibit a $\theta_k$
satisfying~\eqref{nlg_eq:claim2i} and~\eqref{nlg_eq:claim2ii}. 

Consider the $\theta_k$ defined in~(\ref{nlg_eq:6}). By~(\ref{nlg_eq:approx}),
we get~(\ref{nlg_eq:claim2ii}). That $\theta_k$ also
satisfies~(\ref{nlg_eq:claim2i}) is seen using the $L^2$-orthogonal
projection $\mathcal{Q}_k$ into the space of continous functions on
$\Gamma$ that are linear on each $\tau_i^k$. It is well known that
$\mathcal{Q}_k$ is stable in the
$\Hone(\Gamma)$-norm~\cite{BPS02}. Hence, using the standard approximation properties of the
projection $\mathcal{Q}_k\theta$ and the linear interpolant
$\tilde{\theta}_k$ of $\theta$, we find
\begin{align}
\nonumber
 |\tilde\theta_k|_{\Hone(\Gamma)} &\leq |\tilde\theta_k -\mathcal{Q}_k\theta |_{\Hone(\Gamma)} + |\mathcal{Q}_k\theta |_{\Hone(\Gamma)}\\
\nonumber
 &\leq C \, h_k^{-1} \, \|\tilde\theta_k - \mathcal{Q}_k \theta\|_{L^2(\Gamma)}
 + |\mathcal{Q}_k\theta |_{\Hone(\Gamma)}\\
\nonumber
 & \leq C h_k^{-1} \big(  \|\tilde\theta_k-\theta\|_{L^2(\Gamma)} 
    + \| \theta - \mathcal{Q}_k\theta\|_{L^2(\Gamma)} \big)
    + 
\, | \mathcal{Q}_k \theta |_{\Hone(\Gamma)}\\
\label{nlg_eq:9}
& \leq C\, \|\theta \|_{\Hone(\Gamma)}.
\end{align}
%Here we have used the standard approximation properties of the
%projection $\mathcal{Q}_k\theta$ and the linear interpolant
%$\tilde{\theta}_k$ of $\theta$.  
Moreover, by the Sobolev
inequality~(\ref{nlg_eq:Sobolev}), we also have
\begin{equation}
  \label{nlg_eq:7}
  \|\theta_k\|^2 = \sum_i |\theta(x_i)|^2 \|\phi_i\|^2 \leq
  \|\theta\|_{\Hone(\Gamma)}^2 \,\left(\sum_i\int_{\tau_i}1 \right)
  \leq C\,\|\theta\|_{\Hone(\Gamma)}^2.
\end{equation}
Therefore, combining~(\ref{nlg_eq:9}) and (\ref{nlg_eq:7}), we have
\begin{align*}
\| \theta_k \|_{A_k}^2 
& = \| \theta_k\|_{L^2(\Gamma)}^2 + | \tilde{\theta}_k
|_{\Hone(\Gamma)}^2  
 \le C \, \| \theta\|_{\Hone(\Gamma)}^2, 
\end{align*}
which proves~(\ref{nlg_eq:claim2i}).

To complete the proof of the first inequality in \eqref{nlg_eq:innerpr} of the lemma,
\begin{align*}
\|\sigma\|_{\Hmone(\Gamma)}^2 & = \sup_{\theta\in\Hone(\Gamma)}\frac{(\sigma , \theta)^2}{\|\theta\|_{\Hone(\Gamma)}^2} =\sup_{\theta\in\Hone(\Gamma)}\frac{(\sigma , \varPi_k\theta)^2}{\|\theta\|_{\Hone(\Gamma)}^2}\\
&\leq C \sup_{\theta\in\Hone(\Gamma)}\frac{(\sigma , \varPi_k\theta)^2}{(A_k \varPi_k \theta ,\varPi_k\theta)}=C\, (A_k^{-1} \sigma ,\sigma),
\end{align*}
where we have used~(\ref{nlg_eq:QkProj}).
}

%%%%%%%%%%%%%%%%%%%%%%%%%%%%%%%%%%%%%%%%%%%%%%%%%%%%%%%%%%%%%%
%%    smoothers
%%%%%%%%%%%%%%%%%%%%%%%%%%%%%%%%%%%%%%%%%%%%%%%%%%%%%%%%%%%%%%
\subsection{Smoothers}\label{nlg_ssec:smooth}

Using the inner products from Section~\ref{nlg_sec:discr} we can
define a simple Richardson smoother suitable for multigrid algorithms
(see for example \cite{B93}). The smoother is given by
\begin{align}
\label{nlg_eq:4}
[R_k \sigma,\theta]_{k}&=\frac{1}{\tilde{\lambda}_k}(\sigma,\theta)_{-1},
\end{align}
where $\tilde{\lambda}_k$ is  related to  the Rayleigh-Ritz
quotient involving the {\em positive definite} sesquilinear form
$\Lambda(\cdot\, ,\, \cdot)$, namely
\begin{align}\label{nlg_EVP}
\lambda_k&=\sup_{\theta\in\mathcal{M}_k} \frac{\Lambda(\theta, \theta)}{[\theta,\theta]_{k}},
\end{align}
as follows:  We assume for our analysis  that
$\tilde{\lambda}_k$ is a number that satisfies
\begin{equation}
\label{nlg_eq:20}
\lambda_k \le \tilde{\lambda}_k \le C \lambda_k  
\end{equation}
for some mesh independent constant $C$.

Define the operator $\Lambda_k$ as in~(\ref{nlg_eq:1}) but with the
sesquilinear form $\mathcal{V}(\cdot ,\cdot)$ replaced by
$\Lambda(\cdot ,\cdot)$.  The eigenvalue $\lambda_k$ is then a
computable version of the largest eigenvalue of the operator
$\Lambda_k$ with respect to the minus one inner product. In practice,
we could choose $\tilde{\lambda}_k = \lambda_k$, or an approximation
to $\lambda_k$ computed by a few iterations of the power method.

%%%%%%%%%%%%%%%%%%%%%%%%%%%%%%%%%%%%%%%%%%%%%%%%%%%%%%%%%%%%%%
%%    Matrix version 
%%%%%%%%%%%%%%%%%%%%%%%%%%%%%%%%%%%%%%%%%%%%%%%%%%%%%%%%%%%%%%
\subsection{Matrix version}\label{nlg_ssec:matrix}

Now we give a readily implementable matrix version of the previously
given multigrid algorithm. 
Since $\mathcal{M}_k \subseteq \mathcal{M}_{k+1}$ we can find numbers
$c_{i,l}$ such that $\phi^k_i=\sum_{l=1}^{N_{k+1}}\,c_{i,l}\,
\phi^{k+1}_l$. These entries define the $N_k \times N_{k+1}$
restriction matrix $\mathtt{C}_k$ by $[\mathtt{C}_k]_{i,l} = c_{i,l}$.
This matrix and its transpose are used as intergrid transfer operators.
We further define the operator
\begin{align}
\label{nlg_eq:3}
  \mathtt{f}_k
  &: \mathcal{M}_k\longrightarrow \mathbb{C}^{N_k}, 
  && [\mathtt{f}_k (\sigma)]_i=  (\sigma\, ,\,\phi^k_i)_{-1}.
\end{align}
Algorithm \ref{nlg_a:1} can then be translated into an approximation
scheme for the matrix version $ \mathtt{V_J} \mathtt{u} = \mathtt{b}$
of equation \eqref{nlg_eq:2}. Here, the vectors are given by $\mathtt{b} =
\mathtt{f}_k(f_J)$ and $\mathtt{u} = \mathtt{e}_k(\sigma_J)$ and the
system matrix is $\mathtt{V_J}=[\langle V_J \phi^J_j , \phi^J_i
  \rangle]_{i,j}$.  This yields a procedure
$\mathtt{Mg}_J(\mathtt{s},\mathtt{b})$ that outputs an approximation
to the solution given an input iterate $\mathtt{s}$.  To describe it
we will also need a matrix of the operator $A_k$, which we denote by
$\mA_k$.  Specifically $\mA_k$ is the matrix satisfying
$\cf_k(A_k\sigma)=\mA_k \, \cf_k(\sigma)$. It is a circulant matrix
with cyclically shifted rows of the form
\[
\mA_k=\mathop{\mathrm{Circulant}} 
\begin{bmatrix}
0\cdots 0,&
-(l_i\,l_{i-})^{-1}, &
1+l_i^{-2}+(l_i\,l_{i-})^{-1},&-l_i^{-2},
&0\cdots 0 
\end{bmatrix}.
\]
Note that $\mA_k$ is neither tridiagonal nor symmetric, but $\mH_k
\mA_k$ is symmetric, where $\mathtt{H}_k$ is a diagonal matrix whose
$i^{\mathrm{th}}$ diagonal entry is  $\mathrm{meas}( \tau^k_i )$.
The translation of Algorithm~\ref{nlg_a:1} into its matrix version is done via
the identities of the following lemma. 
\begin{lemma}
  \label{nlg_lem:matrix}
  The following identities hold:
\begin{align}
\label{nlg_Nigam_contrib_eq:7}
  \ev_k\,(V_k\, g)
   & =\mV_k \,\cf_k( g),
  && \text{ for all } g\in\mathcal{M}_k,
  \\
\label{nlg_Nigam_contrib_eq:8}
  \ev_{k-1}({Q_{k-1}}\, g )
  & =\mC_{k-1}\,\ev_k(g)
  && \text{ for all } g\in\mathcal{M}_{k},
  \\
\label{nlg_Nigam_contrib_eq:9}
  \cf_k(g )
  & =\mC^t_{k-1} \cf_{k-1}(g )
  && \text{ for all } g\in\mathcal{M}_{k-1},
  \\
\label{nlg_Nigam_contrib_Nigam_contrib_eq:10}
  \cf_k(R_k\,g )
  & =\tilde{\lambda}_k^{-1} \,
  \mH_k^{-1}\,\mA_k^t\, \ev_k(g)
  && \text{ for all } g\in\mathcal{M}_{k},
  \\
\label{nlg_Nigam_contrib_Nigam_contrib_eq:11}
  \cf_1(V_1^{-1}g)
  & = \mV_1^{-1}\ev_1(g)
  && \text{ for all } g\in\mathcal{M}_{1}.
\end{align}
\end{lemma}
\proof{
  Let us prove~(\ref{nlg_Nigam_contrib_eq:7}):
   \begin{align*}
     [\ev_k\,(V_k\, g)]_i 
     & = (V_kg, \phi^k_i)_{-1} =  \mathcal{V}(g, \phi^k_i)
     \\
     & = \sum_{j=1}^{N_k}  [\cf_k(g)]_j 
     \mathcal{V}( \phi^k_j, \phi^k_i) 
     = \sum_{j=1}^{N_k}
     [\mV_k]_{i,j}   [\cf_k(g)]_j  \\
     &= [\mV_k \cf_k(g)]_i.
   \end{align*}
   Next, let us prove~(\ref{nlg_Nigam_contrib_eq:8}):
 \begin{align*}
   [\ev_{k-1}({Q_{k-1}}\, g )]_i 
   & = (Q_{k-1} g, \phi^{k-1}_i )_{-1} 
    =   (g, \phi^{k-1}_i )_{-1} 
    =   (g, \sum_{l=1}^{N_k} [\mC_{k-1}]_{i,l} \phi^k_l )_{-1} 
   \\
   & =   [\mC_{k-1} \ev_k(g)]_i,
 \end{align*}
 since $\mC_{k-1}$ is real. The proof of~(\ref{nlg_Nigam_contrib_eq:9}) is
 similar. To prove~(\ref{nlg_Nigam_contrib_Nigam_contrib_eq:10}), observe
 that
\begin{align*}
  \frac{1}{\tilde{\lambda}_k} [\ev_k(g)]_i \equiv
  \frac{1}{\tilde{\lambda}_k}(g,\phi^k_i)_{-1}
  = 
  [R_k g,\phi^k_i]_{k} & = [\mH_k \mA_k^{-1}\,  \cf_k(R_k g)]_i.
\end{align*}
Multiplying both sides by the symmetric matrix $\mA_k \mH_k^{-1}$, we obtain~(\ref{nlg_Nigam_contrib_Nigam_contrib_eq:10}).
Proofs of the other identities are similar.  
}

These identities enable us to state a matrix version of
Algorithm~\ref{nlg_a:1}. For example, applying $\cf_k$ to the
step~(\ref{nlg_eq:13}) of Algorithm~\ref{nlg_a:1} and using
\eqref{nlg_Nigam_contrib_Nigam_contrib_eq:10} and using
Lemma~\ref{nlg_lem:matrix}, we have
\begin{eqnarray*}
\cf_k(\sigma_1)		 
& =&\cf_k(\sigma)	 + \cf_k(R_k(f-V_k\sigma) )  
 = \cf_k(\sigma)	 + 
\tilde{\lambda}_k^{-1} \,
  \mH_k^{-1}\,\mA_k^t\, \ev_k(f-V_k\sigma) )  \\
& = &\cf_k(\sigma)	 + 
\tilde{\lambda}_k^{-1} \,
  \mH_k^{-1}\,\mA_k^t\, (\ev_k(f)- \mV_k\cf_k(\sigma) ).
\end{eqnarray*}
Thus, the matrix version of this step is $ \mathtt{s}_1 = \mathtt{s} +
\tilde{\lambda}_k^{-1} \, \mH_k^{-1}\,\mA_k^t\, (\mathtt{b} -
\mV_k\mathtt{s}) $ with $\mathtt{s}_1 = \cf_k(\sigma_1), \mathtt{s}=
\cf_k(\sigma), $ and $\mathtt{b} = \ev_k(f)$.  Using also the other
identities
in~(\ref{nlg_Nigam_contrib_eq:7})--(\ref{nlg_Nigam_contrib_Nigam_contrib_eq:11}),
we can similarly translate the entire algorithm. We then obtain the
following matrix version of the algorithm
$\Mgm_J(\mathtt{s},\mathtt{b})$, which outputs an
approximation for the solution of the matrix equation $ \mV_J
\mathtt{u} = \mathtt{b}$, given an input iterate $\mathtt{s}$.
\begin{algorithm}
  \label{nlg_A:2}
  Let $\mathtt{s}$ and $\mathtt{b}$ be any given vectors in $\mathbb{C}^{N_k}$. Define
  $\Mgm_k(\mathtt{s},\mathtt{b})$ recursively as follows.
  Set $\Mgm_1(\mathtt{s},\mathtt{b})=\mV_1^{-1}\mathtt{b}$. If $k>1$, define
  $\Mgm_{k}(\mathtt{u}, \mathtt{b})$ as the vector in  $\mathbb{C}^{N_k}$ obtained
  recursively by:
  \begin{align*}
    \mathtt{s}_1 & = \mathtt{s} +
    \tilde{\lambda}_k^{-1} \, \mH_k^{-1}\,\mA_k^t\, (\mathtt{b} - \mV_k\mathtt{s}),
    \\
    \Mgm_k(\mathtt{s},\mathtt{b})& =\mathtt{s_1} + \mC_{k-1}^t\Mgm_{k-1}(\mathtt{0}, 
    \mC_{k-1} (\mathtt{b}- \mV_k\mathtt{s_1})). 
  \end{align*}
\end{algorithm}
It is important to note that the inverse of $\mA_k$ is not needed in
the implementation. We only need to multiply by $\mA_k^t$.  In the
case of a uniform mesh with mesh size $h_k$ the multiplication by the
matrix $\mH_k^{-1}$ reduces to multiplication with the constant
$1/h_k$. Furthermore, it has been shown in \cite{BLP94} that
\begin{equation}
\label{nlg_eq:18}
\lambda_k=O(1/{h_k}).
\end{equation}
Motivated by this observation, we could replace the expression
$\tilde{\lambda}_k^{-1}\,\mH_k^{-1}$ in Algorithm~\ref{nlg_A:2} by $C$,
for some constant $C$ (which depends on the geometry of the
domain). This would allow us to bypass the eigenvalue computation,
which would otherwise be inherent in the algorithm. The numerical
experiments we present later, however, explicitly include the
expression $\tilde{\lambda}_k^{-1}\,\mH_k^{-1}$.

Note that a matrix preconditioner $\mB_k$ for $\mV_k$ is implicit in
Algorithm~\ref{nlg_A:2} and is defined by $\mB_k \mathtt{b} =
\Mgm_k(0,\mathtt{b})$.  A theoretical study of the convergence rate of
the algorithm is presented in the next section.

%%%%%%%%%%%%%%%%%%%%%%%%%%%%%%%%%%%%%%%%%%%%%%%%%%%%%%%%%%%%%%
%%    convergence analysis    
%%%%%%%%%%%%%%%%%%%%%%%%%%%%%%%%%%%%%%%%%%%%%%%%%%%%%%%%%%%%%%
\section{Convergence Analysis}\label{nlg_sec:conv}

%%%%%%%%%%%%%%%%%%%%%%%%%%%%%%%%%%%%%%%%%%%%%%%%%%%%%%%%%%%%%%
%%    preliminary steps    
%%%%%%%%%%%%%%%%%%%%%%%%%%%%%%%%%%%%%%%%%%%%%%%%%%%%%%%%%%%%%%

\subsection{Preliminary Steps}

Before we can give a detailed description of the convergence behavior
of Algorithm \ref{nlg_a:1}, we need to pave the way with some preliminary
discussions.  In particular, we need the Galerkin projections $P_k:
\Hmhalf(\Gamma)\longrightarrow\Mk$  satisfying
\begin{equation}
\label{nlg_eq:Pk}
\mathcal{V}(P_k \sigma,\mu)=\mathcal{V}(\sigma,\mu) \qquad\mbox{ for
  all } \mu\in\Mk.  
\end{equation}
As we see next, such operators are well defined once the mesh size is
sufficiently small. The assumption on the wave number that $\kappa^2$
is not an interior eigenvalue of $-\Delta$ implies that for homogeneous
right hand side the equation $V\sigma=0$ only has the trivial solution
$\sigma=0$.  It is then a standard theorem on compact perturbations
(see for example \cite[Theorem 4.2.9]{SS04}), that the discretized
version of equation \eqref{nlg_eq:sl} has a unique solution
$\sigma_k\in\Mk$ if the corresponding meshsize $h_k$ is sufficiently
small.  Furthermore, we know that the Galerkin solutions $\sigma_k$
converge quasi-optimally to the true solution $\sigma$, i.e.
\[
\|\sigma -\sigma_k\|_{\Hmhalf(\Gamma)}\;\leq\;
C\;\min_{v_k\in\Mk}\|\sigma-v_k\|_{\Hmhalf(\Gamma)}.
\]
Also, once the mesh size is sufficiently small, the Galerkin solutions
depend continuously on the data, i.e.,
\begin{align}\label{nlg_E:cont}
 \|\sigma_k\|_{\Hmhalf(\Gamma)}\;\leq\; C \;\|f \|_{\Hhalf(\Gamma)}.
\end{align}
As a consequence, we immediately have the following lemma which shows
that $P_1$ is a well defined continuous operator once $h_1$ is small
enough (and so is $P_k$ for $k>1$). In the lemma and elsewhere, we
write $\|\cdot\|_{\Lambda}$ for the vector norm ${\Lambda(\cdot\,
  ,\cdot)}^{1/2}$. We will also use the same notation for the operator
norm induced by  this vector norm.

\begin{lemma}   \label{nlg_lem:Pk}
  There exists an $H>0$ such that once the coarse mesh size $h_1$ is
  less than $H$, there is a unique $P_k\sigma$
  satisfying~\eqref{nlg_eq:Pk} for all $k\ge 1$ and moreover,
  \begin{align}
    \|P_k\sigma\|_{\Lambda} \leq C\;\|\sigma\|_{\Lambda}.
  \end{align}
\end{lemma}

Now let us introduce a few ingredients needed to analyze the multigrid
algorithm.  A simple induction argument shows that $\Mg_k(\cdot\,
,\cdot)$ as defined in Algorithm~\ref{nlg_a:1} is linear as a mapping from
$\Mk\times\Mk$ into $\Mk$. It is also consistent in the sense that
$\sigma_k=\Mg_k(\sigma_k,V_k\sigma_k)$ holds for all
$\sigma_k\in\Mk$. The error reduction operator of the scheme is given
by \begin{align} \E=\Mg_J(\cdot\, , 0),
\end{align}
i.e., if $e^i$ denotes the error at step $i$, we have
$e^{i+1}=\Mg_J(e^i,0)$.  Furthermore, the error reduction operator
admits a product representation as shown in
Lemma~\ref{nlg_lemma:prod}. This representation will be essential in the
convergence analysis of the V-cycle scheme. Proofs of such results can
be found in~\cite{B93}.

\begin{lemma}\label{nlg_lemma:prod}
Let $\Tk=\Rk\Vk\Pk$ for $k\geq 2$ and set $T_1=P_1$. For $k\geq 1$ we then define $E_k u=u-\Mg_k(0,\Vk\Pk u)$ and set $E_0=I$, the identity operator. Then, 
\begin{align}
E_k&=E_{k-1}(I-\Tk),\qquad\text{ and }\\
\E \;&=(I-T_1)(I-T_2)\cdots(I-T_J) . \label{nlg_E:prod}
\end{align}
\end{lemma}

The same representation holds for the error reduction operator
$\Etilde$ of the definite problem.  Analogous to~(\ref{nlg_eq:Pk}), we can
define $\tilde{P}_k$ as the orthogonal projection into $\Mk$ with
respect to the $\Lambda(\cdot,\cdot)$-inner product.  This is the
Galerkin projection for the principal part of the differential
operator.  If we set $\tilde{T}_k=\Rk \Lambda_k \Pktilde$, we get
\begin{align}\label{nlg_E:prod2}
\Etilde=(I-\tilde{T}_1)\cdots(I-\tilde{T}_J).
\end{align}
This operator is proved to be a reducer in~\cite{BLP94}.
Specifically, in~\cite{BLP94} the convergence for the symmetric
version of the multigrid algorithm applied to the positive definite
problem was shown.  In fact, it was shown that the symmetric error
reduction operator $\tilde{E}^s$ in this case is bounded away
from $1$ independently of the number of levels of refinement.  The
symmetric version differs from Algorithm~\ref{nlg_a:1} by an additional
post-smoothing step. However, it is well known (see,
e.g.,~\cite[Remark 3.4]{BP93}) that the analogous result holds for
Algorithm \ref{nlg_a:1} with just the pre-smoothing, i.e., we have the
following theorem.
 \begin{theorem}\label{nlg_thm:posdef}
The error reduction operator $\Etilde$ for Algorithm~\ref{nlg_a:1} applied to the positive definite problem satisfies
\begin{align}
\|\Etilde\|_{\Lambda}& \leq \tilde\delta < 1,
\end{align} % $\eqref{nlg_E:prod2}$
where $\tilde\delta$ is independent of $J$.
 \end{theorem}

In order to analyze the algorithm for the indefinite Helmholtz case we
look at the difference between $\E$ and $\Etilde$. Let
$Z_k=T_k-\tilde{T}_k$, and suppose for some positive $\alpha$ we have
\begin{align}
\|Z_k\|_{\Lambda}	&\leq C_1\,h_k^{\alpha}		&\text{for } k=1,\ldots , J. \label{nlg_eq:zk}
\end{align}
With this assumption, by well known arguments in an abstract multigrid
setting~\cite{B93,GPD04}, we have the following theorem.
\begin{theorem}\label{nlg_thm:pert}
Let $\E$ satisfy \eqref{nlg_E:prod} and $\Etilde$ satisfy
\eqref{nlg_E:prod2}. Assume that \eqref{nlg_eq:zk} holds.  Then,
there exists a positive constant $C_2$ depending on $C_1$, $h_1$, and
$\alpha$ above, such that:
\begin{align}\label{nlg_E:errindef}
\|\E\|_{\Lambda}\leq ||\Etilde||_{\Lambda}+C_2 \, h_1^{\alpha}.
\end{align}
\end{theorem}

We know that $\| \Etilde \|_\Lambda \le \tilde\delta<1$ by
Theorem~\ref{nlg_thm:posdef}. Hence by virtue of Theorem~\ref{nlg_thm:pert},
to prove a convergence result for our multigrid application, we only
need to verify the hypotheses of Theorem~\ref{nlg_thm:pert},
namely~(\ref{nlg_eq:zk}). This will be done in
Subsection~\ref{nlg_ssec:convergence}.

Before concluding this subsection, we need to establish one more
ingredient for the multigrid perturbation argument.  It is well known
that the difference between the single layer potentials of the
Helmholtz and the Laplace equations, namely $D:=V-\Lambda$, is compact
as a map $\Hmhalf(\Gamma)\rightarrow\Hhalf(\Gamma).$ For our purposes we need the following lemma.

\begin{lemma} \label{nlg_lemma:mappingD}
$D$ is bounded as a map $\Hmhalf(\Gamma)\longrightarrow\Hone(\Gamma)$.
\end{lemma}
\proof{In this proof, we use the explicit integral representation of
  the single layer potentials as given in Section~\ref{nlg_intro}. The operator $D =
  V - \Lambda$ generates the sesquilinear form
  \[
  D( \mu,\sigma) = \langle D \mu, \sigma \rangle,
  \]
  and  is an integral operator whose kernel consists of
  the function
  \[
  f(x,y) = g(|x-y|)
  \quad\text{ where } \quad
  g(z) = 
  \frac{i}{4}H_0^{(1)}(\kappa \, z) +\frac{1}{2\pi}
  \ln(z).
  \]
  The function $g$ has the following asymptotic behavior as $z$
  approaches $0$:
  \begin{align}
    \label{nlg_eq:10}
    g(z)  & \sim  c_1 + O(z^2 \log z) \\
    \label{nlg_eq:11}
    g'(z) & \sim  c_2 (z \log z) + O( z) \\
    \label{nlg_eq:15}
    g''(z) & \sim c_3 + O(\log z ),
  \end{align}
  for some constants $c_i$ (depending on $\kappa$).

  Let us now estimate the $H^1$-norm of $D\sigma$. Denote by $\partial
  f$ the derivative of $f$ with respect to arc length along
  $\Gamma$. Then
\begin{align}
 \|D \sigma \|_{\Hone(\Gamma)}^2 &= \|D \sigma \|_{L^2(\Gamma)}^2 \; +\;
 \|\partial( D\sigma)\|_{L^2(\Gamma)}^2.
 \label{nlg_terms}
\end{align}
Let us start by bounding the first term on the right hand side.
Letting $F_x(y) \equiv f(x,y)$, we have 
\begin{align}
 \|D \sigma \|_{L^2(\Gamma)}^2	&= \int_{\Gamma}\;\left| \int_{\Gamma}
 f(x,y)\sigma(y)\; ds_y \;\right|^2  ds_x\notag\\	
\label{nlg_aux1}
 &\leq \int_{\Gamma} 
 \| F_x \|^2_{\Hhalf(\Gamma)}\|\sigma\|^2_{\Hmhalf(\Gamma)}\;  ds_x.
\end{align}
By the trace theorem
\[
\| F_x \|_{\Hhalf(\Gamma)} \le C \| F_x \|_{\Hone(B_r)},
\]
for some ball $B_r$ of sufficiently large radius $r$ (so that $B_r$
contains $\Gamma$). The term $\| F_x \|_{\Hone(B_r)}$ is finite because
\begin{equation}
  \label{nlg_eq:17}
\nabla_y F_x = g'(|x-y|) \; \frac{x - y }{|x-y|},
\end{equation}
is a square integrable function (due to~(\ref{nlg_eq:11}) and the
boundedness of $(x - y )/|x-y|$).

By a change of variable (mapping $x$ to $0$), integrals of $F_x$ can
be converted to integrals of $F_0$ on transformed domains. Hence, by
enlarging the transformed integration region, we have
\[
\| F_x \|_{\Hone(B_r)} \le \| F_0 \|_{\Hone(B_{2r})}.
\]
This shows that the first factor in the integrand of~(\ref{nlg_aux1})
admits a bound independent of the integration variable~$x$, so
\begin{align}
  \nonumber
  \|D \sigma \|_{L^2(\Gamma)}^2	
  & \le \int_\Gamma \| F_x
  \|^2_{\Hhalf(\Gamma)} \|\sigma\|^2_{\Hmhalf(\Gamma)}\;  ds_x
   \\
\nonumber
  &\le \meas(\Gamma) \, C \,\| F_0 \|_{\Hone(B_{2r})} 
  \|\sigma\|^2_{\Hmhalf(\Gamma)}
  \\ \label{nlg_eq:D0}
  & \le
  C  \,\|\sigma\|^2_{\Hmhalf(\Gamma)}.
\end{align}
  
We treat the second term in \eqref{nlg_terms} similarly. Denote by $t_x$ the unit tangential vector to $\Gamma$ in the point $x\in\Gamma$, which is defined everywhere except on a set of measure zero (the corners). Then, 
\begin{align}
\nonumber
 \|\partial D \sigma \|_{L^2(\Gamma)}^2	&= \int_{\Gamma}\;\left|
 \left( \nabla_x \int_{\Gamma} f(x,y)\sigma(y) \; ds_y \right)\cdot
 t_x  \;\right|^2  ds_x\\
\nonumber
		& =	\int_{\Gamma}\;\left|  \int_{\Gamma}
\left(\nabla_x f(x,y) \right) \cdot t_x \; \sigma(y) \; ds_y
\;\right|^2  ds_x\\
\label{nlg_aux2}
&\leq \int_{\Gamma} \| G_x \|^2_{\Hhalf(\Gamma)}
\|\sigma\|^2_{\Hmhalf(\Gamma)}\;  ds_x
\end{align}
where
\[
G_x(y) = t_x \cdot \nabla_x f(x,y) = -t_x \cdot \nabla_y F_x.
\]
Note that by differentiating~(\ref{nlg_eq:17}),
\[
\nabla_y G_x = -g''( |x-y|) \frac{ (x-y)\cdot t_x}{|x-y|^2} (x-y)
\,-\,
g'(|x-y|) \nabla_y\bigg(\frac{ (x-y)\cdot t_x }{ |x-y|}\bigg).
\]
The term $\nabla_y(x-y)\cdot t_x / |x-y|$ is $O(1/|x-y|)$, while the
multiplying factor $g'(|x-y|)$ is $ O( |x-y| \log|x-y|)$
by~\eqref{nlg_eq:11}. Hence the last term is $O(\log|x-y|)$. The first
term on the right hand side is also $O(\log|x-y|)$, because
of~\eqref{nlg_eq:15}. Consequently, $\nabla_y G_x$ is locally square
integrable on $\R^2$. Therefore, returning to~\eqref{nlg_aux2}, we can
complete the estimation using a trace inequality and bounding $\| G_x \|^2_{\Hhalf(\Gamma)}$ 
independently of $x$ as before. Thus
\begin{equation}
  \label{nlg_eq:16}
 \|\partial D \sigma \|_{L^2(\Gamma)}^2 \le C
 \|\sigma\|^2_{\Hmhalf(\Gamma)}.
\end{equation}
Using~(\ref{nlg_eq:16}) and~\ref{nlg_eq:D0} in~(\ref{nlg_terms}), the proof is finished.
}

%%%%%%%%%%%%%%%%%%%%%%%%%%%%%%%%%%%%%%%%%%%%%%%%%%%%%%%%%%%%%%
%%    Convergence   
%%%%%%%%%%%%%%%%%%%%%%%%%%%%%%%%%%%%%%%%%%%%%%%%%%%%%%%%%%%%%%
\subsection{Convergence}    \label{nlg_ssec:convergence}

Now we give our main result on the convergence of the multigrid
algorithm for our application. The proof proceeds by verifying the
hypotheses of Theorem~\ref{nlg_thm:pert}. For this, we need a regularity
result. Consider the solution $\veps$ of the adjoint problem 
\begin{equation}
  \label{nlg_eq:eps}
  \VV( \eta, \veps ) = \overline{F(\eta)} \qquad\text{ for all } \eta
  \in H^{-1/2}(\Gamma),
\end{equation}
for some linear functional $F$ on $ H^{-1/2}(\Gamma)$, or in other
words $F$ is in $\Hhalf(\Gamma)$. If $F$ is more regular, then we
expect the solution $\veps$ to be more regular. 

To make this precise, note that~\eqref{nlg_eq:eps} can be rewritten as
\[
\VVb( \veps,\eta ) = F(\eta),
\]
where  $\VVb(\cdot,\cdot)$ is defined for smooth $\sigma,\mu$ by
\[
\VVb(\sigma,\mu) = \int_{\Gamma}\!\int_{\Gamma} 
\overline{
\frac{i}{4}
\Hnkl(\kappa |x-y|)}\, \sigma (y) \,\overline{\mu(x)} \; ds_y  \; ds_x.
\]
This form extends continuously to $H^{-1/2}(\Gamma) \times
H^{-1/2}(\Gamma)$ and the operator $V^*: \Hmhalf(\Gamma) \mapsto
\Hhalf(\Gamma)$ defined by $\ip{V^* \sigma,\mu} = \VVb(\sigma,\mu)$ is
continuous. It can be written as
\[
V^* = \Lambda + D^*
\]
where $D^*$ is an integral operator analogously to $D$, but with an
integral kernel conjugate to that of $D$. The same type of arguments as in
Lemma~\ref{nlg_lemma:mappingD} show that
\[
D^* : \Hmhalf(\Gamma) \mapsto H^1(\Gamma)
\]
is continuous.
Now, it is well known~\cite[Thm. 3.2.2]{SS04} that %\cite{cs89}
for the positive definite problem
$\Lambda u = F$, there is a  regularity result:
\[
\| u \|_{H^s(\Gamma)} \le C \| F \|_{H^{s+1}(\Gamma)}
\qquad\text{ for }  0\le s < s_0
\]
where $s_0$ is a positive number depending only on the angles of
$\Gamma$. Applying this result with $s=0$ to~\eqref{nlg_eq:eps} rewritten
as $V^*\veps = F$, or in other words, $\Lambda \veps = F - D^* \veps$,  we obtain that
\[
\begin{aligned}
\| \veps \|_{L^2(\Gamma)} 
& \le C ( \| F \|_{\Hone(\Gamma)} + \| D^* \veps \|_{\Hone(\Gamma)} )\\
& \le C ( \| F \|_{\Hone(\Gamma)} + \| \veps \|_{\Hmhalf(\Gamma)} )
\end{aligned}
\]
by the above mentioned continuity of $D^*$.  Now by the unique
solvability of~\eqref{nlg_eq:eps}, we also have the stability estimate
$\|\veps \|_{\Hmhalf(\Gamma)} \le C \| F \|_{\Hhalf(\Gamma)}.$ This,
together with the continuous imbedding of $\Hone(\Gamma)$ into
$\Hhalf(\Gamma)$ shows that
\begin{equation}
  \label{nlg_eq:reg}
\| \veps \|_{L^2(\Gamma)} \le C \| F \|_{\Hone(\Gamma)}.
\end{equation}
We will use this regularity result in the proof of the next theorem.

\begin{theorem}\label{nlg_thm:convergence}
There is an $H>0$ and a  $0<\delta <1$ such that whenever the coarse
grid meshsize $h_1$ is less than $H$, the error reduction operator $E$
of Algorithm~\ref{nlg_a:1} applied to the indefinite acoustic single layer
equation satisfies
\begin{align}
\|\E \|_{\Lambda}& \leq \delta.
\end{align}
Here, $\delta$ is independent of the refinement level $J$.
\end{theorem}
\proof{ 
This proof proceeds by verifying~\eqref{nlg_eq:zk} and applying
  Theorem~\ref{nlg_thm:pert}.  To verify \eqref{nlg_eq:zk}, we begin with the
  following:
\begin{align}
| \mathcal{D}(\sigma ,\mu) | & = |\langle D \sigma ,  \mu\rangle |
 \leq \|D\sigma\|_{\Hone(\Gamma)} \|\mu\|_{\Hmone (\Gamma)}\notag\\
&\leq\; C \;\|\sigma\|_{\Hmhalf(\Gamma)}\|\mu\|_{\Hmone(\Gamma)}.\label{nlg_eq:19}
\end{align}
This is a consequence of Lemma~\ref{nlg_lemma:mappingD}. We shall
use~\eqref{nlg_eq:19} several times below.

We first prove~\eqref{nlg_eq:zk} for $k>1$. Define $\Dktilde = \Vk\Pk -
\Lambdak\Pktilde$. Then 
$\langle \Dktilde \sigma,\mu_k\rangle =
\mathcal{D}(\sigma,\mu_k)$ for all $\sigma $ in $\Hmhalf(\Gamma)$ and all
$\mu_k$ in $M_k$. 
\[
\Zk=\Tk-\Tktilde=\Rk\left( \Vk\Pk -
\Lambdak\Pktilde\right)=\Rk\Dktilde.
\]
For any $\sigma\in\mathcal{M}_J$ and $k>1$, we have
\begin{align*}
 \|Z_k \sigma\|_{\Lambda}^2 & = 
\Lambda (\Rk\Dktilde\sigma , \Rk\Dktilde\sigma)\leq
\;\tilde{\lambda}_k\;[\Rk\Dktilde\sigma , \Rk\Dktilde \sigma]_k 
&&\text{ by~\eqref{nlg_eq:20}} \\
&=\;\tilde{\lambda}_k\;\frac{1}{\tilde{\lambda}_k} (\Dktilde\sigma ,
\Rk\Dktilde\sigma )_{-1}=\mathcal{D}(\sigma , \Rk\Dktilde\sigma)
&&\text{ by~\eqref{nlg_eq:4}}\\
&\leq \; C\; \|\sigma\|_{\Hmhalf(\Gamma)}\|\Rk\Dktilde\sigma\|_{\Hmone(\Gamma)}
&&\text{ by~\eqref{nlg_eq:19}}.
\end{align*}
The last factor  can be estimated by
\begin{align*}
 \|\Rk\Dktilde\sigma\|_{\Hmone(\Gamma)}^2&=(\Rk\Dktilde\sigma ,\Rk\Dktilde\sigma )_{-1}\leq \; C \; [\Rk\Dktilde\sigma, \Rk\Dktilde\sigma]_k\\
				&= \; C\; \frac{1}{\tilde{\lambda}_k} (\Dktilde\sigma ,\Rk\Dktilde\sigma )_{-1}\leq \; C\; \frac{1}{\tilde{\lambda}_k} \|\Dktilde\sigma\|_{\Hmone(\Gamma)}\|{\Rk\Dktilde\sigma}\|_{\Hmone(\Gamma)},
\end{align*}
and also noting that
\begin{align*}
 \|\Dktilde\sigma\|_{-1}^2 	&= (\Dktilde\sigma ,\Dktilde\sigma)_{-1}= \left((\Vk\Pk-\Lambdak\Pktilde)\sigma ,\Dktilde\sigma\right)_{-1}\\
				&=\mathcal{D}(\sigma ,\Dktilde\sigma)\leq \; C\; \|\sigma\|_{\Hmhalf(\Gamma)}\|\Dktilde\sigma\|_{\Hmone(\Gamma)}.
\end{align*}
In combination, these show 
\begin{align*}
 \|\Zk\sigma\|_{\Lambda}=\|\Rk\Dktilde\sigma\|_{\Lambda}\leq \; C\; \tilde\lambda_k^{-1/2}\;\|\sigma\|_{\Hmhalf(\Gamma)}\;\leq\;C\; h_k^{{1/2}}\;\|\sigma\|_{\Hmhalf(\Gamma)}.
\end{align*}
The last inequality follows from~\eqref{nlg_eq:20} and~\eqref{nlg_eq:18}.
Hence, for $k\geq 2,$ 
\begin{align*}
 \|\Zk\|&=\sup_{\sigma\in\Mk}\;\frac{\|
   \Zk\sigma\|_{\Lambda}}{\|\sigma\|_{\Lambda}}\;\leq\; C\;
 h_k^{{1/2}},
\end{align*}
so we have verified \eqref{nlg_eq:zk} with $\alpha=1/2$.

To prove~\eqref{nlg_eq:zk} on the coarsest level ($k=1$), we will
use~\eqref{nlg_eq:reg} and the following duality argument along the lines
of a similar argument in~\cite{CS88}. Let $\sigma$ in $\Mk$. Define
\[
F(\eta) = (\eta, \,(I-P_1)\sigma)_{-1}.
\]
This is a continuous linear functional on $\Hmhalf(\Gamma)$ and hence
there is a unique solution $\veps$ to~\eqref{nlg_eq:eps} with this
$F$. Hence,
\begin{align}
\nonumber
\| (I-P_1)\sigma \|_{\Hmone(\Gamma)}^2 
& =  \overline{ F( (I-P_1)\sigma ) }\\
\nonumber
& = \VV(  (I-P_1)\sigma, \veps )\\
\nonumber
& = \mathcal V ( (I - P_1)\sigma, \veps - \veps_1  ) \\\label{nlg_eq:23}
& \le C \| (I - P_1)\sigma\|_{\Hmhalf(\Gamma)} \|  \veps - \veps_1\|_{\Hmhalf(\Gamma)}
\end{align}
for any $\veps_1$ in $\mathcal{M}_1$. We choose an $\veps_1$ with
optimal approximation properties. Note that since
\[
\| F \|_{H^1(\Gamma)} = \| (I - P_1)\sigma \|_{\Hmone(\Gamma)},
\]
the regularity result~\eqref{nlg_eq:reg} holds for $\veps$. Therefore,
\begin{align*}
\| \veps - \veps_1\|_{\Hmhalf(\Gamma)} & \le C h_1^{1/2}\| \veps
\|_{L^2(\Gamma)}\le C h_1^{1/2} \| F \|_{\Hone(\Gamma)}\\
& = C h_1^{1/2} \| (I - P_1)\sigma \|_{\Hmone(\Gamma)}.
\end{align*}
Using this in~\eqref{nlg_eq:23}, we conclude that 
\begin{equation}
\label{nlg_eq:21}
  \| (I - P_1)\sigma\|_{\Hmone(\Gamma)} \le C h_1^{1/2}   \| (I - P_1)\sigma\|_{\Hmhalf(\Gamma)} .
\end{equation}

We use~\eqref{nlg_eq:21} to estimate the norm of $Z_1$ as follows.
\begin{align*}
\Lambda(Z_1 \sigma, \mu_1)
&=
\Lambda((P_1-\tilde{P}_1)\sigma\, , \, \mu_1)
=D\left((I-P_1)\sigma\, ,\, \mu_1\right)\\
& \leq C\, \|(I-P_1)\sigma\|_{\Hmone(\Gamma)}\, \|\mu_1\|_{\Lambda} 
&&\text{ by~\eqref{nlg_eq:19}}\\
& \leq C\, h_1^{1/2}\|(I-P_1)\sigma\|_{\Hmhalf(\Gamma)}\, \|\mu_1\|_{\Lambda} 
&&\text{ by~\eqref{nlg_eq:21}}\\
& \leq C\, h_1^{1/2}\|\sigma\|_{\Hmhalf(\Gamma)}\, \|\mu_1\|_{\Lambda} 
&&\text{ by~Lemma~\ref{nlg_lem:Pk}}.
\end{align*}
This proves~\eqref{nlg_eq:zk} for the $k=1$ case as well. 

Hence, we can apply Theorem~\ref{nlg_thm:pert} with $\alpha = 1/2$ to get 
\[
\| \E \|_\Lambda  \le \| \Etilde \|_\Lambda + C_2 h_1^{1/2} \le
\tilde\delta + C_2 h_1^{1/2},
\]
where $\tilde\delta$ is a positive number less than one given by
Theorem~\ref{nlg_thm:posdef}.  It is now clear that when $h_1$ is small
enough, the result follows.}

%%%%%%%%%%%%%%%%%%%%%%%%%%%%%%%%%%%%%%%%%%%%%%%%%%%%%%%%%%%%%%
%%     numerical experiments     
%%%%%%%%%%%%%%%%%%%%%%%%%%%%%%%%%%%%%%%%%%%%%%%%%%%%%%%%%%%%%%
\section{Numerical Experiments}
\label{nlg_sec::ne}

\subsection{Effect of weaker inner product on eigenfunctions}

The performance of the multigrid algorithm described in Section~\ref{nlg_sec::ma} depends crucially on the spectral behavior of the positive definite operator $\Lambda$. For a discretized version of this integral operator the eigenfunctions corresponding to small magnitude eigenvalues are highly oscillatory, while those eigenfunctions corresponding to the large end of the spectrum are non-oscillatory. 
Standard multigrid approaches are successful for operator equations with the opposite spectral behavior and the use of the weaker inner products from 
Section~\ref{nlg_sec:discr} effectively transforms the single layer problem into this setting.
In this section we present the details of two examples showing the undesirable behaviour of the discrete eigenfunctions of the stiffness matrix associated with $\Lambda$ while working in the (natural) $\Hmhalf(\Gamma)$ inner product, and the effect of working with the weaker inner product instead. We demonstrate the effectiveness of this approach for two geometries, a (smooth) circle and a Lipschitz domain (square).

Recall again that if the boundary $\Gamma$ is discretized by means of a partition $ x_1, x_2,$ $\ldots ,x_N ,x_{N+1}=x_1$, we denote by $\phi_i$ and $l_i$ the characteristic function and respectively the length of the element $\tau_i=\mathrm{conv}(x_i,x_{i+1})$. The span of the $\left\{\phi_i\right\}$ is denoted by $\mathcal{M}_k$. 
In the case of the square, the boundary discretization consists of straight line segments; whereas in the case of the circle our discretization consists of arcs of equal angle.
The basis coefficients of an element   $\sigma\in\mathcal{M}_k$ with respect to the $\left\{\phi_i\right\}$ are given in terms of the vector $\cf(\sigma)$ defined in \eqref{nlg_eq:5}.

\begin{figure}
\begin{center}\includegraphics[width=0.7\textwidth]{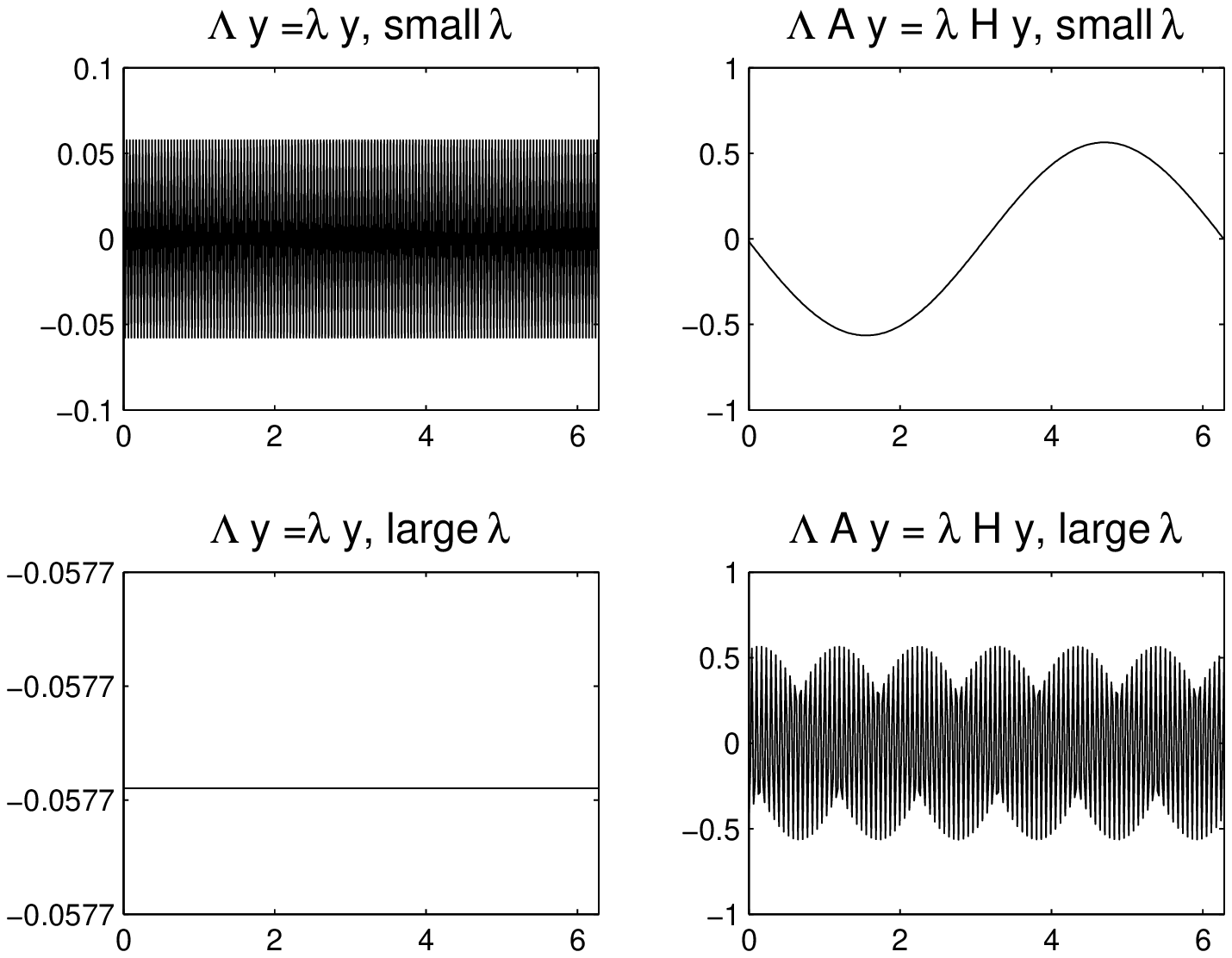}\end{center}
\caption{Eigenfunction behavior, circular curve}{Eigenfunctions for a circular curve with \mbox{$N=300$}. They correspond to \\
the smallest e-val of \mbox{$\mLa$} (top left), the smallest e-val of \mbox{$\mLa\mA y =\lambda \mH y$} (top right),
the largest e-val of \mbox{$\mLa$} (bottom left) and the largest e-val of \mbox{$\mLa\mA y =\lambda \mH y $} (bottom right). }
\label{nlg_fig:Circle}
\end{figure}

For a given discretization we are interested in the spectrum of the $N\times N$ stiffness matrix corresponding to the single layer operator for the Laplacian with entries $[\mLa]_{i,j}=\langle \Lambda \phi_j\, ,\, \phi_i\rangle$.
In Figure~\ref{nlg_fig:Circle} the left-hand plots show the eigenfunction of $\mLa$ for a circular domain, corresponding to the smallest (top left) and largest (bottom left) eigenvalues respectively. These figures illustrate the phenomenon described above. 

In Section~\ref{nlg_sec:discr} we introduced a discrete inner product on the discretization space and proved that its associated norm is equivalent to the natural $\Hmone(\Gamma)$ norm. As before we denote by $\mA$ the finite difference matrix corresponding to $-u''+u$ on $\Gamma$ with periodic boundary conditions, and by $\mH$ the diagonal matrix with $i^{\mathrm{th}}$ diagonal entry $l_i$.
According to the definition of the discrete inner product in \eqref{nlg_eq:12} we find for $\sigma,\theta\in\mathcal{M}_k$:
\begin{align*}
 [\sigma ,\theta]_k 	&=  \langle \sum_i \,[\cf(A^{-1}\sigma)]_i\,\phi_i \;,\;\sum_j \,[\cf(\theta)]_j\,\phi_j\rangle_{\Gamma}
			 = \sum_i\; [\cf(A^{-1}\sigma)]_i \;\overline{([\cf(\theta)]_i)}\; l_i\\
			&= [\cf(\theta)]^* \,\mH\,\mA^{-1}\,[\cf(\sigma)],
\end{align*}
and also,	
\begin{align*}
\Lambda (\theta,\theta)&=\sum_{i,j}[\cf(\theta)]_i \;\overline{([\cf(\theta)]_j)}\;\underbrace{\Lambda(\phi_i,\phi_j)}_{\mLa_{j,i}}=[\cf(\theta)]^*\;\mLa\;[\cf(\theta)]. 
\end{align*}
The smoothing procedure, defined in Section~\ref{nlg_ssec:smooth}, depends on the largest eigenvalue of the following Rayleigh quotient with respect to the weaker inner product.
\begin{align*}
 \lambda &=\sup_{\theta\in\mathcal{M}_k}\;\frac{\Lambda(\theta,\theta)}{[\theta,\theta]_k} =\sup_{\theta\in\mathcal{M}_k}\;\frac{\cf(\theta)^*\;\mLa\;\cf(\theta)}{\cf(\theta)^* \,\mH\,\mA^{-1}\,\cf(\theta)}=\sup_{y\in\mathbb{R}^N}\;\frac{y^*\;\mA^* \mLa\;\mA \;y}{y^* \mA^*\,\mH\,y}.
\end{align*}
The two matrices $\mA^*\,\mLa\,\mA$ and $\mA^*\,\mH$ are Hermitian. Therefore, $\lambda$ is the largest generalized eigenvalue of the problem:
\begin{align}
 \mA^*\,\mLa\,\mA \, y	&=  \lambda \;\mA^*\,\mH \, y,\qquad \mbox{ or equivalently, }\qquad
\mLa\,\mA \, y			=  \lambda \; \mH \, y.  \label{nlg_GenEV}
%\mLa\,\mA \, \mH^{-1}\,z 	&=  \lambda \; z 	&\mathrm{with} \;z=\mH y \notag\\
\end{align}
The corresponding eigenfunction is given in terms of its basis coefficients $\cf(\theta) =\mA\,y$.
We note that the $L^2$ norm of $\theta$ is easy to compute.
\begin{align*}
 \|\theta\|_{L^2}^2 = \sum_{i,j} [\cf(\theta)]_i \overline{[\cf(\theta)]_j}\; \langle \phi_i ,\phi_j\rangle = \cf(\theta)^*\,\mH\,\cf(\theta) = y^* \,\mA^*\,\mH \,\mA \,y .
\end{align*}
This allows us to normalize the coefficient vector $\cf(\theta)$ such that $\theta$ has norm one in $L^2(\Gamma)$. 
We are now ready to examine the spectral behavior in the smoothing operation in terms of the generalized eigenvalue problem \eqref{nlg_GenEV}. 

In the case where $\Gamma$ is a circle centered at the origin, the entries of the matrix $\mLa$ are particularly simple to compute. We require the radius to be bounded by $R<1/2$ in order to guarantee the positive definiteness of the integral operator $\Lambda$ and we discretize the circle with arcs $\Gamma_i$, $i=1..N$ of equal angle. 

\begin{align}
[\mLa]_{i j}  & = \langle\Lambda \Phi_j ,\Phi_i \rangle = -\frac{1}{2\pi} \int_{\Gamma_i} \int_{\Gamma_j} \log(|x-y|) ds_x ds_y \notag\\
	& = -\frac{1}{2\pi} \,  \frac{1}{2} \int_{\Gamma_i} \int_{\Gamma_j} \log(2 R^2 \left(1-\cos(\theta_i -\theta_j) \right)) \,R^2\, d\theta_j d\theta_i  \notag\\
	& = -\frac{R^2}{4\pi}\, \int_{\Gamma_i} \int_{\Gamma_j} \log(2 R^2) \, d\theta_j d\theta_i  \;-\; \frac{R^2}{4\pi}\underbrace{\int_{\Gamma_i} \int_{\Gamma_j} \log \left(1-\cos(\theta_i -\theta_j) \right)\, d\theta_j d\theta_i }_{I_{ji}} \notag\\
	& =  -\frac{R^2}{4\pi}\log(2 R^2) \, \left(\frac{2\pi}{N}\right)^2 \; -\; \frac{R^2}{4\pi}\, I_{ji}  = -\frac{R^2}{4\pi}\,\left\{ \log(2 R^2) \, \left(\frac{2\pi}{N}\right)^2\;+\; I_{ji}\right\}
\end{align}
Hence, we have to evaluate the weakly singular integrals \[I_{ji}=\int_{\Gamma_i} \int_{\Gamma_j} \log \left(1-\cos(\theta_i -\theta_j) \right))\, d\theta_j d\theta_i\] 
for all possible choices of $\Gamma_i$ and $\Gamma_j$. Due to symmetry of the circle, the single layer matrix is a symmetric Toeplitz matrix and hence it is sufficient to compute its first row. For $|i-j|>1$, the integrand in $I_{ij}$ is smooth, and we use Gaussian quadrature to compute the entries. 
In the two remaining cases when either $|i-j| = 1 \mod N$ or $i=j$, we separate the singularity of the integrand according to 
\begin{align}
\log(1-\cos(t)) & = \log(\frac{t^2}{2}-\frac{t^4}{4!}+\frac{t^6}{6!}-\frac{t^8}{8!}+\ldots)\notag\\
		& = \underbrace{\log(t^2)}_{f_s} \;  \underbrace{-\; \log(2) \; +\; \log(1 -\frac{2t^2}{4!}+\frac{2t^4}{6!}+\ldots)}_{f_a}\label{nlg_eqn:singularity}\\
		& = f_s(t) \; +\; f_a(t).\notag
\end{align}
We then integrate the singular term $f_s(t)$ exactly and use Gaussian quadrature to compute the integral over the nonsingular function $f_a(t)$.

Figure~\ref{nlg_fig:Circle} shows the effect of the smoothing operator. The eigenfunction of the generalized eigenvalue problem $\mLa \mA y =\lambda H y$ corresponding to the smallest eigenvalue (top right) is smooth (and represents a function of period $2\pi$ on the circle. The eigenfunction of the largest eigenvalue (bottom right) is highly oscillatory. All eigenfunctions shown correspond to $N=300$ elements on the circle.

\begin{figure}
\begin{center}\epsfig{file=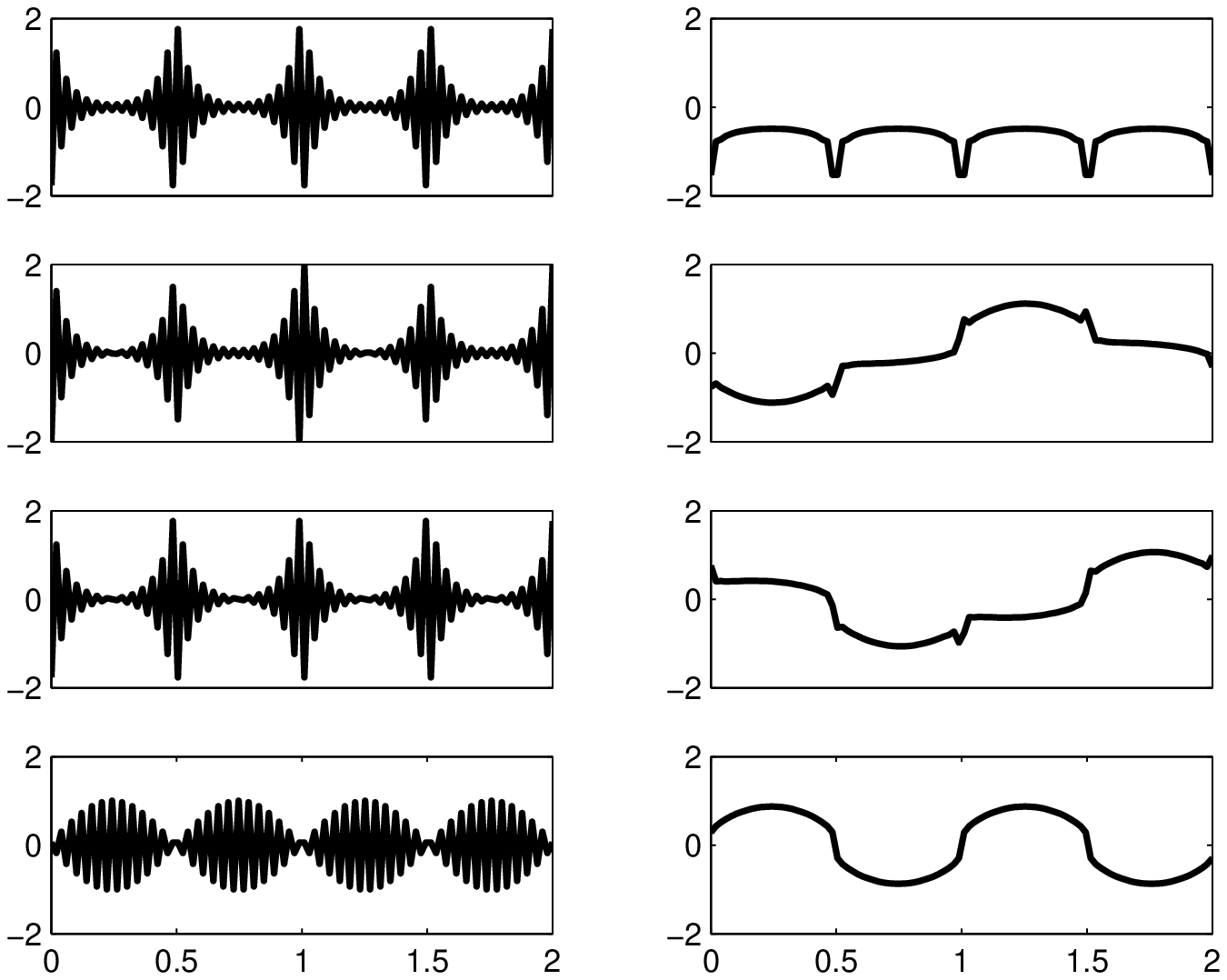, width=0.7\textwidth}\end{center}
\caption{Eigenfunction behavior, square boundary curve, small eigenvalues}{Eigenfunction behavior for a square boundary curve $\Gamma$ and a uniform mesh with mesh size \mbox{$h=1/50$}. Eigenfunctions correspond to the four smallest \\
e-vals of $\mLa$ (left) and the four smallest generalized e-vals of \mbox{$\mLa\mA y =\lambda \mH y $} (right).}
\label{nlg_fig:SquareSmall}
\end{figure}

\begin{figure}
\begin{center}\epsfig{file=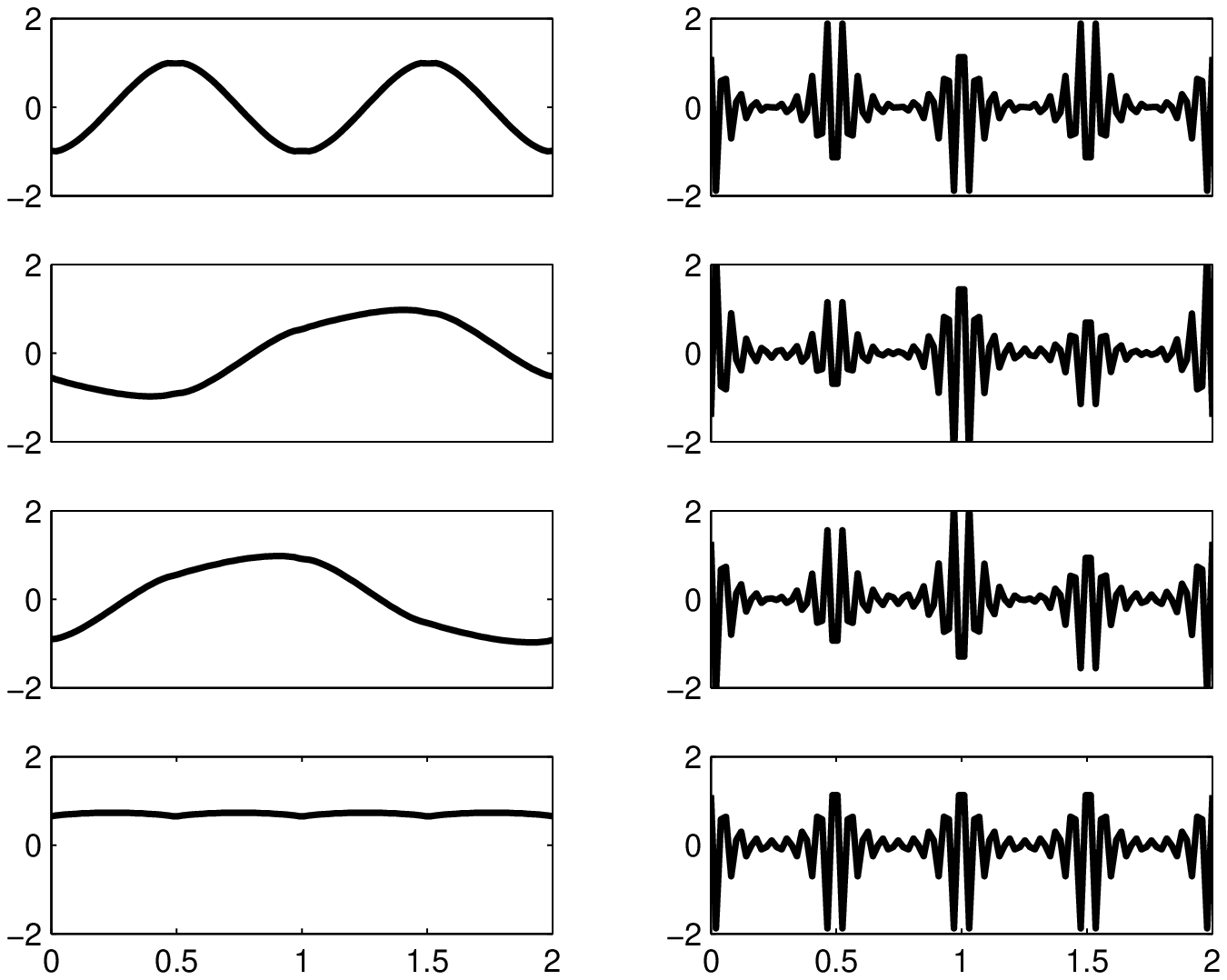, width=0.7\textwidth}\end{center}
\caption{Eigenfunction behavior, square boundary curve, large eigenvalues}{Eigenfunction behavior for a square boundary curve $\Gamma$ and a uniform mesh with mesh size \mbox{$h=1/50$}. Eigenfunctions correspond to the four largest \\
e-vals of $\mLa$ (left) and the four largest generalized e-vals of \mbox{$\mLa\mA y =\lambda \mH y $} (right).}
\label{nlg_fig:SquareLarge}
\end{figure}

As a second test case we consider a square with side length $1/2$. 
We compute the eigenfunctions corresponding to the four smallest and the four largest eigenvalues for both the stiffness matrix $\mLa$ and the generalized eigenvalue problem $\mLa\mA y =\lambda \mH y$. Figures~\ref{nlg_fig:SquareSmall}
and \ref{nlg_fig:SquareLarge} illustrate the results for a uniform discretization of the boundary courve $\Gamma$ with meshsize $h=1/50$. 
Again, we observe that the eigenfunctions of the generalized eigenvalue problem display the reversed (and sought-after) smoothness behavior.
The same behavior has been observed in cases of quasi-uniform meshes.

% Figures~\ref{nlg_fig:SquareSmallNonUni}
% and \ref{nlg_fig:SquareLargeNonUni} correspond to a semi-uniform discretization of the boundary courve $\Gamma$ with meshsize $h_1=1/100$ on two opposite sides of the square and $h_2=1/120$ on the two other sides.
%It is worth noting that our computations we Further computations have been done for example for semi-uniform discretizations of the boundary courve %$\Gamma$ with meshsize $h_1$ on two opposite sides of the square and $h_2$ on the two other sides. 

% \begin{figure}
% \begin{center}\epsfig{file=sq_4small_nonuniform.eps, width=0.7\textwidth}\end{center}
% \caption{Eigenfunction behavior for a square boundary curve $\Gamma$ with a semi-uniform mesh of meshsize $h=1/100$ on two opposite sides of the square and $h=1/120$ on the two other sides.
%  Left: eigenfunctions corresponding to the four smallest eigenvalues of $\mLa$. 
% Right: eigenfunctions corresponding to the four smallest generalized eigenvalues which satisfy $\mLa\mA y =\lambda \mH y $.}
% \label{nlg_fig:SquareSmallNonUni}
% \end{figure}
% 
% 
% \begin{figure}
% \begin{center}\epsfig{file=sq_4large_nonuniform.eps, width=0.7\textwidth}\end{center}
% \caption{Eigenfunction behavior for a square boundary curve $\Gamma$ with a semi-uniform mesh of meshsize $h=1/100$ on two opposite sides of the square and $h=1/120$ on the two other sides. Left: eigenfunctions corresponding to the four largest eigenvalues of $\mLa$. Right: eigenfunctions corresponding to the four largest generalized eigenvalues which satisfy $\mLa\mA y =\lambda \mH y $.}
% \label{nlg_fig:SquareLargeNonUni}
% \end{figure}

\subsection{Multigrid convergence results}
In this section we present numerical convergence results for the multigrid algorithm described in Section~\ref{nlg_sec::ma} to underline its effective use. Various results from a different set of experiments have already been reported in \cite{GGN08}.
Recall that we want to solve the exterior Helmholtz problem with prescribed Dirichlet data on the boundary of a scattering object. The fundamental solution $i/4\;\Hnkl(\kappa\|x-x^*\|)$ solves the Helmholtz equation away from its singularity $x^*$ and satisfies the correct growth condition at infinity. If we place the singularity into the interior of the scattering domain, we can use this so called point source problem as a convenient test case, for which the exact solution is known. We also present results for the scattering of an incident plane wave by polygonal shaped obstacles. In all the tables below, {\bf H} is the coarsest mesh size, and {\bf h} is the finest mesh.

The fast implementation of the matrix-vector multiplications in the algorithm was not subject of our study and hence we do not report on the overall CPU time used by the algorithm. We anticipate that the CPU time will be competitive once those matrix operations are implemented using matrix-compression techniques such as H-matrices. 

\subsubsection{Effect of domain shape on performance}
We first present results for point-source scattering from 4 different objects: a square, a rectangle (sides of ratio 1:4), an equilateral triangle and a thin wedge.  The thin wedge is described in terms of the xy-coordinates of its three corner points, namely $(0,0)$, $(1/3, \sqrt{3}/3)$ and $(0,1/15)$. This amount to an angle of $\pi/3$ between the horizontal x-axis and the lower edge of the wedge. In each of these examples, the diameter of the object is less than 1, which guarantees the positive definiteness of the potential single layer operator in the sense of \eqref{nlg_eq:posdef}. In each of Tables \ref{TablesObject1}-\ref{TablesObject4}, the proposed multigrid scheme is used as a linear solver. We report the number of iteration numbers required to reach a given relative residual norm. The point source is located inside the domain, so the true solution is known in each case.

\begin{table}
\caption{Linear multigrid iteration counts with $\kappa=2.1$,  point-source inside square,$10^{-6}$ relative residual norm}
\label{TablesObject1}
\begin{tabular}{cc|c|c|c|c|c|c}
\multicolumn{7}	{c|}{ \bf		 H	}		 		& 	Degrees	\\
	& &1/2		&1/4		&1/8		&1/16	&1/32 	& of freedom\\
	\hline
	&1/4		&17	& -	& -	& -	&- & 32\\
	&1/8		&16	& 15	& -	&-	& -& 64\\
	&1/16	& 15& 15 & 15 & -	& -& 128\\
{\bf h}&1/32 	&15 & 15 & 15 & 15 & -&256\\
	& 1/64	&16 & 16 & 15 & 15 & 15& 512\\
	& 1/128	& 16&16 &16 &16 	&16 &1024\\
	& 1/256	& 16&16 &16 &16 	&16 &2048
\end{tabular}
\end{table}

\begin{table}
\caption{Linear multigrid iteration counts with $\kappa=2.1$, point-source inside  rectangle,$10^{-6}$ relative residual norm}
\label{TableObject2}
\begin{tabular}{cc|c|c|c|c|c|c}
\multicolumn{7}	{c|}{ \bf		 H	}		 & Degrees	\\
	& 		&1/2	&1/4	&1/8	&1/16&1/32 	& of freedom\\
	\hline
	&1/4		&15	& -	& -	& -	&- 		& 72\\
	&1/8		&17	& 15	& -	&-	& -		& 144\\
	&1/16	& 18& 18 & 16 & -	& -		& 288\\
{\bf h}&1/32 	&18 & 19 & 18 & 16 & -		&576\\
	& 1/64	&19 & 19 & 19 & 18 &16 		&1152\\
	& 1/128	&19 &19 &19 & 19	& 19		&2304
\end{tabular}
\end{table}
\begin{table}
\caption{Linear multigrid iteration counts with $\kappa=2.1$, point-source inside triangle,$10^{-6}$ relative residual norm}
\label{TableObject3}
\begin{tabular}{cc|c|c|c|c|c|c|c|c|c}
\multicolumn{10}	{c|}{ \bf		 H	}		 	& Degr. of\\
	& 		&1&1/2	&1/4	&1/8	&1/16&1/32 &1/64&1/128	& freedom\\
	\hline
	&1/2		&21& -	& -	& -	&- 	&-	&-&-		& 6\\
	&1/4    	&33&30&-&-&-&-&-&-				&12\\
	&1/8		& 37& 37	& 36	&-	& - &-&-&-			& 24\\
	&1/16	& 28& 28 & 27 & 27	& -	&-	&-&-		& 48\\
{\bf h}&1/32 	& 24& 24 & 24 & 23 & 23&-	&-&-		&96\\
	& 1/64	& 23 & 23 & 23 & 22 &22 &22	&-&-		&192\\
	& 1/128	& 22 &22 & 22 & 22	& 22	&21	&22&-	&384\\
	& 1/256	&21&21&21&21&21&21&21&21		&768\\
	& 1/512 	&20&20&20&21&21&21&20&20		&1536
\end{tabular}
\end{table}

\begin{table}
\caption{Linear multigrid iteration counts with $\kappa=2\pi$, thin wedge, pre-smoothing only,
$10^{-4}$ relative residual norm}
\label{TablesObject4}
\begin{tabular}{cc|c|c|c|c|c|c|c|c|c}
\multicolumn{7}	{c|}{ \bf		 H	}		 	& Degr. of\\
	& 		&1/4& 1/8	& 1/16	& 1/32	& 1/64 & freedom\\
	\hline
	&1/16   	& 31&21&-&-&-				&352\\
{\bf h}&1/32	&18 &13 &9 & - &	-	& 704\\
	& 1/64	& 17& 11 &8  &5 		&-		& 1408\\
	&1/128 	& 17& 11 & 7 &4   &3			&2816\\
\end{tabular}
\end{table}

\subsubsection{Multigrid as linear solver/ preconditioner}
In this section, we provide results about the use of the proposed multigrid scheme as a linear solver (Tables \ref{TableMG1,TableMG2}), and  as a preconditioner for GMRES (used without restart, Tables \ref{TableGMRES1,TableGMRES2}). We present iteration counts in each case with just presmoothing, or with both pre-and post-smoothing.  We use the same equilateral triangle as the domain as in the previous section. Again, a point-source is placed inside the domain, so we can compare with the exact solution. In contrast to the previous section, here we present results for wave number $\kappa=10.2$. The tolerances in relative residual norm are $10^{-6}, 10^{-9}$ for the linear solver and preconditioned GMRES, respectively. For both, multigrid as a linear solver and preconditioned GMRES, the iteration numbers stay almost constant with increasing degress of freedom, whereas the iteration numbers of GMRES grow considerably. Convergence for the algorithm with an additional postsmoothing step follows from our convergence result by standard arguments in multigrid theory.

\begin{table}
\caption{Linear multigrid iteration counts with $\kappa=10.2$, point source inside triangle,
pre-smoothing only,
$10^{-6}$ relative residual norm}
\label{TableMG1}
\begin{tabular}{cc|c|c|c|c|c|c|c}
\multicolumn{8}	{c|}{ \bf		 H	}		 	& Degr. of	\\
	& 		&1/4	&1/8	&1/16&1/32 &1/64&1/128	& freedom\\
	\hline
	&1/8		& 38& -	& -	&-	& - &			& 24\\
	&1/16	& 36& 33 & - & -	& -	&-		& 48\\
{\bf h}&1/32 	& 32& 30 & 28 & - & -&-			&96\\
	& 1/64	& 27& 26 & 25 & 24 &- &-			&192\\
	& 1/128	& 24 &23& 22 & 22	& 22	&-		&384\\
	& 1/256	&22&21&21&21&21&22			&768\\
	& 1/512 	&21&21&21&21&21&21			&1536
\end{tabular}
\end{table}
\begin{table}
\caption{Linear multigrid iteration counts with $\kappa=10.2$,point source inside  triangle, 
pre- and post-smoothing,
$10^{-6}$ relative residual norm}
\label{TableMG2}
\begin{tabular}{cc|c|c|c|c|c|c|c}
\multicolumn{8}	{c|}{ \bf		 H	}		 	& Degr. of	\\
	& 		&1/4	&1/8	&1/16&1/32 &1/64&1/128	& freedom\\
	\hline
	&1/8		& 19& -	& -	&-	& - &			& 24\\
	&1/16	& 20& 17& - & -	& -	&-			& 48\\
{\bf h}&1/32 	& 18&16 &14  & - & -&-			&96\\
	& 1/64	& 17& 15& 14 & 13 &- &-			&192\\
	& 1/128	& 15&14& 13 & 13&13 &-			&384\\
	& 1/256	&15&13&13&13&13&13			&768\\
	& 1/512 	&16&14&13&13&13&13			&1536
\end{tabular}
\end{table}

\begin{table}
\caption{GMRES iteration counts, $\kappa=10.2$,point source inside  triangle,  triangle, pre-smoothing only,
$10^{-9}$ relative residual norm}
\label{TableGMRES1}
\begin{tabular}{cc|c|c|c|c|c}
\multicolumn{6}	{c|}{ \bf		 H	}		 	& GMRES	\\
	& 		&1/2	&1/4	&1/8	&1/16 &without preconditioning\\
	\hline
	&1/8		&17&16&-	& - 			&16\\
	&1/16	&23&21&20&-			&24\\
{\bf h}&1/32 	&23&22&22&19		&31\\
	& 1/64	&24&23&23&22		&37\\
	& 1/128	&25&25&24&23		&44\\
	& 1/256	&25&27&24&23		&52\\
	& 1/512 	&26&28&24&24			&63\\
	& 1/1024	&27&28&26&26			&74
\end{tabular}
\end{table}

\begin{table}
\caption{GMRES iteration counts, $\kappa=10.2$, triangle, uniform grid, $\theta=1$, pre and post-smoothing,
$10^{-9}$ relative residual norm}
\label{TableGMRES2}
\begin{tabular}{cc|c|c|c|c|c}
\multicolumn{6}	{c|}{ \bf		 H	}		 	& GMRES	\\
	& 		&1/2	&1/4	&1/8	&1/16 &without preconditioning\\
	\hline
	&1/8		&18&15&-	& - 			&16\\
	&1/16	&20&17&17&-			&242\\
{\bf h}&1/32 	&20&18&18&17		&31\\
	& 1/64	&21&19&19&18		&37\\
	& 1/128	&21&21&19&19		&44\\
	& 1/256	&22&22&21&20		&52\\
	& 1/512 	&22&24&21&21			&63\\
	& 1/1024	&23&24&21&21			&74
\end{tabular}
\end{table}

\subsubsection{Effect of frequency on performance}
Here we present the effect of increasing $\kappa$ on the number of multigrid iterations taken to achieve a given relative residual error. We place a point source inside a square whose diameter is less than 1. We expect the coarsest mesh required should satisfy the constraint $\kappa {\bf H} \approx$ constant. At least in this example, the method performs well even though this constraint was not strictly satisfied. For example, in Table \ref{Tablef1}, $\kappa H=0.252$,  while in Table \ref{Tablef3}, $\kappa H=1.575$, 
 
\begin{table}
\caption{Linear multigrid iteration counts with $\kappa=2.1$,point source inside square, 
pre- and post-smoothing,
$10^{-6}$ relative residual norm}
\label{Tablef1}
\begin{tabular}{cc|c|c|c|c|c|c|c|c|c|}
\multicolumn{9}	{c|}{ \bf		 H	}		 	& Degr. of	\\
	& 		&1/4	&1/8	&1/16&1/32 &1/64&1/128&1/256& freedom\\
{\bf h}&1/32 	& 12& 12& 12 & -&- &	-		&-&256\\
	& 1/64	&13 &13 & 12 &  12& -&	-		&-&512\\
	& 1/128	&13 &13&13  &13 &12 &	-		&-&1024\\
	& 1/256	&13&13&13&13&13&13			&-&2048\\
	& 1/512 	&14 &14&14&14&14&			13&13&4096
\end{tabular}
\end{table}

\begin{table}
\caption{Linear multigrid iteration counts with $\kappa=10.2$,point source inside square, 
pre- and post-smoothing,
$10^{-6}$ relative residual norm. The method did not converge for {\bf H}=1/4.}
\label{Tablef2}
\begin{tabular}{cc|c|c|c|c|c|c|c|c|c|}
\multicolumn{9}	{c|}{ \bf		 H	}		 	& Degr. of	\\
	& 		&1/4	&1/8	&1/16&1/32 &1/64&1/128&1/256& freedom\\
{\bf h}&1/32 	& *& 15& 13 & -&- &	-		&-&256\\
	& 1/64	&* &14 & 12 &  11& -&	-		&-&512\\
	& 1/128	&* &14&12  &11 &11 &	-		&-&1024\\
	& 1/256	&*&13&12&11&11&11			&-&2048\\
	& 1/512 	&* &12&12&12&12&			12&12&4096
\end{tabular}
\end{table}

\begin{table}
\caption{Linear multigrid iteration counts with $\kappa=50.4$,point source inside square, 
pre- and post-smoothing,
$10^{-6}$ relative residual norm. The method did not converge for {\bf H}$>1/32$.}
\label{Tablef3}
\begin{tabular}{cc|c|c|c|c|c|c|c|c|c|}
\multicolumn{9}	{c|}{ \bf		 H	}		 	& Degr. of	\\
	& 		&1/4	&1/8	&1/16&1/32 &1/64&1/128&1/256& freedom\\
{\bf h}&1/32 	&*&*& *&*&-&	-		&-&256\\
	& 1/64	&* &* & * &  *& -&	-		&-&512\\
	& 1/128	&* &*&*  &* &13 &	-		&-&1024\\
	& 1/256	&*&*&*&*&13 &11		&-&2048\\
	& 1/512 	&* &*&*&*&12&			11&11&4096
\end{tabular}
\end{table}

\begin{figure}[htbp] %  figure placement: here, top, bottom, or page
    \centering\epsfig{file=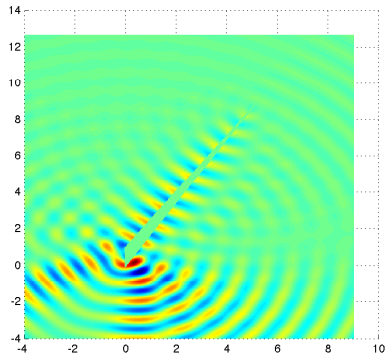, width=0.5\textwidth} 
   \caption{Plane-wave scattering from wedge: linear MG is used to compute the unknown density, and numerical quadrature is used to reconstruct the scattered field.}
   \label{fig:example}
\end{figure}

The numerical experiments show that the multigrid algorithm presented and analyzed in this paper is an efficient tool to solve the first kind single layer equation when used as a linear solver or as a preconditioning procedure for other solvers such as GMRES.


\begin{thebibliography}{10}

\bibitem{BrakhageWerner}
H. Brakhage and P.Werner. 
\newblock Uber das Dirichletsche Aussenraumproblem f\"ur die
Helmholtzsche Schwingungsgleichung. 
\newblock {\em Arch. Math., 16:325Ð329, 1965.}

\bibitem{B93}
J.~H. Bramble.
\newblock {\em Multigrid methods}, volume 294 of {\em Pitman Research Notes in
  Mathematics Series}.
\newblock Longman Scientific \& Technical, 1993.

\bibitem{BKP94}
J.~H. Bramble, D.~Y. Kwak, and Joseph~E. Pasciak.
\newblock Uniform convergence of multigrid {$V$}-cycle iterations for
  indefinite and nonsymmetric problems.
\newblock {\em SIAM J. Numer. Anal.}, 31(6):1746--1763, 1994.

\bibitem{BLP94}
J.~H. Bramble, Z.~Leyk, and J.~E. Pasciak.
\newblock The analysis of multigrid algorithms for pseudodifferential operators
  of order minus one.
\newblock {\em Math. Comp.}, 63(208):461--478, 1994.

\bibitem{BP93}
J.~H. Bramble and J.~E. Pasciak.
\newblock New estimates for multilevel algorithms including the {$V$}-cycle.
\newblock {\em Math. Comp.}, 60(202):447--471, 1993.

\bibitem{BPS02}
J.~H. Bramble, J.~E. Pasciak, and O.~Steinbach.
\newblock On the stability of the {$L\sp 2$} projection in {$H\sp 1(\Omega)$}.
\newblock {\em Math. Comp.}, 71(237):147--156, 2002.

\bibitem{ansatz}
S. Chandler-Wilde and I. Graham.
\newblock Boundary integral methods in high frequency scattering.
\newblock In {\em B. Engquist, A. Fokas, E. Hairer, and A. Iserles, editors, highly oscillatory problems.
Cambridge University Press, 2009}

\bibitem{snc01}
S.N. Chandler-Wilde, I.G. Graham, S. Langdon, and M. Lindner.
\newblock Condition number estimates
for combined potential boundary integral operators in acoustic scattering. Journal of
Integral Equations and Applications, 21:229Ð279, 2009.

\bibitem{snc02} 
S.N. Chandler-Wilde and S. Langdon. 
\newblock A wavenumber independent BEM for an acoustic scattering problem. 
\newblock SIAM J. Numer. Anal., 46:2450Ð2477, 2006.


\bibitem{snc03}
 S.N. Chandler-Wilde and P. Monk.
\newblock Wave-number-explicit bounds in time-harmonic scattering.
\newblock SIAM J. Math. Anal, 39:1428Ð1455, 2008.


\bibitem{CS88}
M.~ Costabel and E.~P. Stephan.
\newblock Duality estimates for the numerical solution of integral equations.
\newblock {\em Numer. Math.}, 54(3):339--353, 1988.

\bibitem{FS97}
S.A. Funken and E.P. Stephan.
\newblock The {BPX} preconditioner for the single layer potential operator.
\newblock {\em Appl. Anal.}, 67(3-4):327--340, 1997.

\bibitem{GGN08}
S.~Gemmrich, J.~Gopalakrishnan, and N.~Nigam.
\newblock The performance of a multigrid algorithm for the acoustic single
  layer equation.
\newblock In {\em Numerical Mathematics and Advanced Applications: Proceedings
  of ENUMATH 2007, the 7th European Conference on Numerical Mathematics and
  Advanced Applications}, pages 175--182. Springer, Heidelberg, 2008.

\bibitem{GPD04}
J.~Gopalakrishnan, J.~E. Pasciak, and L.~F. Demkowicz.
\newblock Analysis of a multigrid algorithm for time harmonic {M}axwell
  equations.
\newblock {\em SIAM J. Numer. Anal.}, 42(1):90--108 (electronic), 2004.

\bibitem{HSW91}
G.C. Hsiao, E.P. Stephan, and W.L. Wendland.
\newblock On the {D}irichlet problem in elasticity for a domain exterior to an
  arc.
\newblock {\em J. Comput. Appl. Math.}, 34(1):1--19, 1991.

\bibitem{HW08}
G.C. Hsiao and W.L. Wendland.
\newblock {\em Boundary integral equations}, volume 164 of {\em Applied
  Mathematical Sciences}.
\newblock Springer-Verlag, Berlin, 2008.

\bibitem{HM73}
George Hsiao and R.~C. MacCamy.
\newblock Solution of boundary value problems by integral equations of the
  first kind.
\newblock {\em SIAM Rev.}, 15:687--705, 1973.

\bibitem{LPR03}
U.~Langer, D.~Pusch, and S.~Reitzinger.
\newblock Efficient preconditioners for boundary element matrices based on
  grey-box algebraic multigrid methods.
\newblock {\em Internat. J. Numer. Methods Engrg.}, 58(13):1937--1953, 2003.

\bibitem{LP07}
U.~Langer and U.~Pusch.
\newblock Convergence analysis of geometrical multigrid methods for solving
  data-sparse boundary element equations.
\newblock In {\em Proc. 8th European Multigrid Conference 2005}. Springer,
  Heidelberg, 2007.

\bibitem{melenk2011}
M. ~L\"ohndorf and J.M.~Melenk 
\newblock Wavenumber-explicit hp-BEM for high frequency scattering
\newblock To appear in {\em SIAM J. Numer. Anal.}, preprint {\tt http://www.asc.tuwien.ac.at/preprint/2010/asc02x2010.pdf}

\bibitem{MMS97}
M.~Maischak, P.~Mund, and E.~P. Stephan.
\newblock Adaptive multilevel {BEM} for acoustic scattering.
\newblock {\em Comput. Methods Appl. Mech. Engrg.}, 150(1-4):351--367, 1997.
\bibitem{melenk2010}
J.M. Melenk and S.Sauter 
\newblock Wavenumber explicit convergence analysis for Galerkin discretizations of the Helmholtz equation
\newblock {\em SIAM J. Numer. Anal. 49 (2011), pp. 1210--1243 }


\bibitem{SS04}
S.~Sauter and C.~Schwab.
\newblock {\em Randelementmethoden, Analyse, Numerik und Implementierung
  schneller Algorithmen}.
\newblock Teubner, 2004.

\bibitem{S03}
O.~Steinbach.
\newblock {\em Numerische N\"aherungsverfahren f\"ur elliptische
  Randwertprobleme}.
\newblock Teubner, Stuttgart, Leipzig, Wiesbaden, 2003.

\end{thebibliography}
\end{document}